\documentclass[twocolumn]{autart}
\usepackage[dvipsnames]{xcolor}
\usepackage{amssymb}
\usepackage{graphicx}
\usepackage{tikz}
\usepackage{mathtools}
\usepackage{arydshln}
\usepackage{cite}
\usepackage[round]{natbib}
\usepackage[colorlinks=true,allcolors=MidnightBlue]{hyperref}
\usepackage{balance}
\allowdisplaybreaks

\edef\endfrontmatter{
    \unexpanded\expandafter{\endfrontmatter}
    \noexpand\endNoHyper 
}

\DeclareMathOperator{\tr}{tr}
\DeclareMathOperator{\sign}{sign}
\DeclareMathOperator{\dist}{dist}
\DeclareMathOperator{\mspan}{span}
\DeclareMathOperator{\blkdiag}{blkdiag}
\DeclareMathOperator{\atan2}{atan2}
\DeclareMathOperator{\re}{Re}
\DeclareMathOperator{\im}{Im}

    \newcommand{\newsharedtheorem}[2]{
    \newtheorem{#1}[mainthm]{#2}
}
\newsharedtheorem{theorem}{Theorem}
\newsharedtheorem{proposition}{Proposition}
\newsharedtheorem{lemma}{Lemma}
\newsharedtheorem{corollary}{Corollary}
\newsharedtheorem{definition}{Definition}
\newsharedtheorem{remark}{Remark}

\def\proof{{\noindent\textbf{Proof.}\ }}
\def\endproof{$\hfill\square$}

\begin{document}

\begin{frontmatter}
\title{An emergence-oriented approach to circular formation\thanksref{footnoteinfo}}

\thanks[footnoteinfo]{This work was partially supported by the Wallenberg AI, Autonomous Systems and Software Program (WASP) funded by the Knut and Alice Wallenberg Foundation.}

\author{Zhaozhan Yao}\ead{zhaozhan@kth.se},
\author{Yuhua Yao}\ead{yuhuay@kth.se},
\author[cor]{Xiaoming Hu}\ead{hu@kth.se}

\corauth[cor]{Corresponding author}
\address{Department of Mathematics, KTH Royal Institute of Technology, Stockholm 100 44, Sweden}
          

\begin{abstract}
In this paper, we study the emergence of circular formation for agents in cyclic pursuit. Each agent is a unicycle traveling at a fixed common forward speed. We first establish a necessary and sufficient condition for the existence of circular formation in cyclic pursuit. Building on this theoretical foundation, we propose a control law that enables the spontaneous formation of circular formations through only local measurements. Notably, key geometric features -- the radius and agent spacing -- are not imposed externally but emerge naturally from the initial conditions of the group. This occurs because the closed-loop system possesses infinitely many non-isolated equilibria, each corresponding to a particular circular formation, and none are asymptotically stable. Consequently, analyzing individual equilibria is no longer informative, and attention is instead directed to the full invariant set (the set of all equilibria). Globally, it is disconnected. Locally, however, each equilibrium together with its neighboring equilibria forms a connected invariant set. This motivates a local stability analysis formulated at the level of invariant sets that are maximally connected. An accompanying stability criterion is then derived and applied to analyze small agent groups ($n \leq 3$), providing insights into the convergence mechanism. Finally, the proposed control law is extended to distance-dependent neighborhoods. Under this setting, the group converges into several clusters, most exhibiting a complete-graph topology. A preliminary stability analysis is then conducted for the case of a complete graph with $n=3$. 
\end{abstract}
\end{frontmatter}

\section{Introduction}
Multi-agent coordination has been an enduring research topic in control community since the beginning of this century.
In particular, owing to the practical potential in applications, circular formation of vehicles with nonholonomic constraints has been actively studied. Theoretical challenges partly stem from the control and analysis of the unicycle model, which is commonly used to describe nonholonomic vehicles.
\\
Early works regarding circular formation are developed in cyclic pursuit and study how simple local rules (using only local measurements) could give rise to complex collective behaviors. In~\cite{marshall2004formations}, unicycles travel at a fixed common forward speed. By employing a steering law where each unicycle's angular speed is proportional to its bearing relative to the next one, the unicycles eventually converge to an equally spaced circular formation. Depending on their spacing (on the circle) and the rotation direction, there are $2(n-1)$ equilibria. Local stability analysis of these equilibria is conducted, revealing which are stable and which are not. A continued study is carried out in~\cite{marshall2006pursuit}, where each unicycle's forward speed is no longer fixed but proportional to its distance from the next one. It is proven that only two special equilibria are locally asymptotically stable, i.e., all unicycles are equally spaced (along the circle) in cyclic order, either clockwise or counterclockwise. This stability, however, holds only for a specific ratio (between the control gains of forward speed and angular speed). Slight deviations from that ratio cause the circular orbit to expand or shrink over time, which makes this behavior non-robust from a practical point of view. To address this issue, a new steering law is proposed in~\cite{zheng2009ring} to guarantee both certain and robust behaviors. The same two equilibria remain stable. Moreover, the ratio no longer needs to be fixed to guarantee stability; rather, it determines the radius of the resulting circle, thus offering greater flexibility. 
\\
Later works shift their focus towards achieving a desired formation, i.e., circular formation control. The (global) prescribed information includes radius/radii, center (stationary or moving), and agent spacing. Some papers are conducted in cyclic pursuit, and the achieved formations are actually variants of equally spaced circular formation~\citep{zheng2015distributed,tripathy2024convergence,galloway2018collective}. The other papers are developed under more general topologies, including complete graphs~\citep{sepulchre2007stabilization,seyboth2014collective}, connected and/or jointly connected graphs~\citep{sun2018circular,sepulchre2008stabilization}, and directed graphs~\citep{yu2022decentralized,sen2019circumnavigation,yu2018cooperative,miao2017cooperative}. These papers consider kinematic unicycles. Extensions to dynamic unicycles have been addressed in~\cite{el2012distributed,shi2021distributed,xie2025circular}. Depending on the problem setup and methodology, these papers establish either local or global stability results. Apart from the unicycle model, some papers instead use single integrators, with the extra assumption that agents can measure bearings~\citep{chun2020multi,zhang2025angle,dou2020target,zhang2020distributed}. 
\\
Typically, the methodologies in most formation control papers take the desired formation as a reference signal and then regulate the formation error. Such an approach, however, encounters limitations in scenarios where the prescribed information is not accessible by certain agents. In recent years, the concept of intrinsic formation control has been applied to produce some special formations like regular polyhedra and antipodal formation~\citep{zhang2020intrinsic,li2022differential}. The point is that no prescribed information of the desired formation is incorporated into the controller; instead, the formation is intrinsically attributed to its geometric properties and the designed topologies. From this perspective, representative cyclic pursuit papers~\cite{marshall2004formations,marshall2006pursuit,zheng2009ring} naturally carry an intrinsic flavor: Identical local control law for all agents, together with the cyclic pursuit scheme, unsurprisingly gives rise to an equally spaced circular formation. \\
Taking a different perspective from prior studies that achieve equally spaced or prescribed circular formation, we view circular formation as an emergent behavior arising from local and self-organized interactions among individuals. A motivating example is fish milling. It has been observed that a fish swarm, when threatened by predators, forms a tightly packed spherical formation, the so-called bait ball. As a simplified problem inspired by this natural phenomenon, in this paper we investigate circular formation from an emergence perspective, exploring how a simple local rule gives rise to complex collective behavior -- spontaneous circular formation. The key geometric features --  radius and agent spacing -- are not fixed or specified a priori, but emerge from the initial conditions of the group. We note that while \cite{sepulchre2007stabilization,sepulchre2008stabilization} realize circular formations with arbitrary agent spacing, the radius remains prescribed and the control laws are based on global rather than local measurements. Overall, the contribution of the paper is threefold:
\\
First, we develop a systematic analysis of  circular formation. A necessary and sufficient condition for the existence of circular formation in cyclic pursuit is derived, revealing the underlying mechanism and informing the control design for the emergence of this type of formations. A set of definitions and algebraic properties are also established. The existence condition readily extends to connected graphs, thereby supporting the subsequent extension to distance-dependent neighborhoods.
\\
Second, a control law using local measurements is proposed and achieves the emergence of circular formation. The closed-loop system possesses infinitely many equilibria, none are isolated and asymptotically stable. Thus, analyzing individual equilibria is no longer meaningful. Instead, attention is directed to the full invariant set (the set of all equilibria), which consists of a finite number of disjoint components (i.e., maximal connected subsets). This motivates a stability analysis formulated at the level of (maximal connected) invariant sets rather than individual equilibria, which stands in clear contrast to most control papers that deal with single or finite isolated equilibria. A local stability criterion is established and a set of stability results are derived for $n\leq3$.
\\
Third, the proposed control law is extended to distance-dependent neighborhoods. For this case, the topology is rendered undirected, and we prove that counterclockwise (resp. clockwise) formations are unstable when the gain parameter is positive (resp. negative). This phenomenon is observed in the simulations of cyclic pursuit as well but cannot be fully analyzed there. Due to the limited sensing range, the group converges to several clusters, most of which have a complete-graph topology. Hence, a preliminary stability analysis is conducted for the case of a complete graph with $n=3$.
\\
The rest of the paper is structured as follows. In section~\ref{sec-problem} we review the relevant equations for unicycles in cyclic pursuit and investigate the mechanism of circular formation. In section~\ref{sec-controller} we propose a control law that achieves spontaneous circular formation and develop a stability analysis approach, together with an accompanying stability criterion. In section~\ref{sec-main} we conclude a set of stability results for $n\leq3$. In section~\ref{sec-global} we slightly modify the control law to achieve almost global convergence for $n=2$. In section~\ref{sec-dist} we extend the control law to distance-dependent neighborhoods, which gives rise to more complex collective behaviors. Finally, section~\ref{sec-conclusion} concludes the paper.
To streamline the presentation, we relegate most of the proofs to the appendices at the end of the paper.

\section{Problem formulation and preliminaries}\label{sec-problem}
In this section, a nonlinear system model for unicycles in cyclic pursuit is reviewed. The target configuration (i.e., circular formation) is defined, and the existence condition and key algebraic properties of circular formation are established, providing a theoretical foundation for the control design and stability analysis.
\subsection{Cyclic pursuit model}\label{sec-model}
Consider a cyclic pursuit of multiple vehicles, where vehicle $i\in\mathcal{I}\coloneqq\{1,2,\ldots,n\}$ pursues $i+1$ modulo $n$. Each is a unicycle with the following nonlinear state model: 
\begin{align}\label{eq-global-kinematic}
    \begin{bmatrix}
        \dot{x}_i \\ \dot{y}_i \\ \dot{\theta}_i
    \end{bmatrix} = \begin{bmatrix}
        \cos{\theta_i} & 0 \\ \sin{\theta_i} & 0 \\ 0 & 1
    \end{bmatrix}\begin{bmatrix}
        v \\ \omega_i
    \end{bmatrix},
\end{align}
where $v>0$ is a fixed forward speed, $z_i= [x_i,y_i]^\top\in\mathbb{R}^2$ denotes the position of unicycle $i$, $\theta_i\in\mathbb{R}$ the heading, and $\omega_i\in\mathbb{R}$ the angular speed to be designed. 
\\
Relative coordinates are considered for analysis convenience. Let $\rho_i$ denote the distance between $i$ and $i+1$, $\alpha_i$ the bearing from $i$'s heading to the heading that would take it directly towards $i+1$, and $\beta_i\coloneqq\theta_i-\theta_{i+1}-\pi$ the heading difference between $i$ and $i+1$ minus $\pi$. After some tedious algebraic manipulations, the motion equations of these coordinates are obtained as:
\begin{subequations}\label{eq-relative-kinematic-1}
    \begin{align}
        \dot{\rho_i} &= -v\left(\cos{\alpha_i}+\cos{(\alpha_i+\beta_i)}\right), \\
        \dot{\alpha_i} &= \frac{v}{\rho_i}\left(\sin{\alpha_i}+\sin{(\alpha_i+\beta_i)}\right)-\omega_i, \\
        \dot{\beta}_i &= \omega_i -\omega_{i+1}.
    \end{align}
\end{subequations}
It is assumed in the relative coordinates that $\rho_i\neq0$. In addition, we exclude any vehicle overlap throughout the paper.
\\
These coordinates can be computed through the following procedure. Let $\bar{z}_i = \left[\bar{x}_i,\bar{y}_i\right]^\top= R(\theta_i)(z_{i+1}-z_i)$, where $R(\theta_i) = \begin{bmatrix}
    \cos{\theta_i} & \sin{\theta_i} \\
    -\sin{\theta_i} & \cos{\theta_i}
\end{bmatrix}$.
Then $\rho_i=\sqrt{\bar{x}_i^2+\bar{y}_i^2}$ and $\alpha_i=\atan2{(\bar{y}_i,\bar{x}_i)}$. Fig.~\ref{fig-rel.coordin.} provides an illustration. Note that $\alpha_i+\beta_i$ is equivalent to the bearing from $i+1$'s heading to the heading that would take it directly towards $i$; therefore, it can be measured in the same manner as $\alpha_i$.
\begin{figure}[h]
    \begin{center}
    \begin{tikzpicture}
    \draw[->, thick] (-.75,0) -- (-.25,0) node[right] {\scalebox{0.75}{$x$}};
    \draw[->, thick] (-.75,0) -- (-.75,.5) node[above] {\scalebox{0.75}{$y$}};
    
    \draw[dashed, thick] (3,.5) -- (4,.5);
    \draw[->, thick] (3,.5) -- ({3+sqrt(2)/2},{.5+sqrt(2)/2});
    \draw[->, dashed, thick] (3.7,.5) arc[start angle=0, end angle=45,radius=0.7];
    \node at (4,.75) {\scalebox{0.75}{$\theta_i$}};
    \node at (2.8,.5) [left] {\scalebox{0.75}{$z_i$}};

    \draw[thick] (3,.5) -- (1,2) node[midway, fill=white, inner sep=1pt] {\scalebox{0.75}{$\rho_i$}};
    \draw[->, thick] ({3+sqrt(2)/2*.5},{.5+sqrt(2)/2*.5}) arc[start angle=45, end angle=142,radius=0.5];
    \node at (2.75,1.2) {\scalebox{0.75}{$\alpha_{i}$}};

    \draw[->, dashed, thick] (3,.5) -- (3,1.5);
    \draw[->, dashed, thick] (3,.5) -- ({3-sqrt(2)/2},{.5-sqrt(2)/2});
    \draw[->, thick] (3,.8) arc[start angle=90, end angle=-135,radius=0.3];
    \node at (3.5,.25) {\scalebox{0.75}{$\beta_{i}$}};

    \draw[dashed, thick] (1,2) -- (2,2);
    \draw[->, thick] (1,2) -- (1,3);
    \draw[->, dashed, thick] (1.5,2) arc[start angle=0, end angle=90,radius=0.5];
    \draw[->, thick] (1,2.3) arc[start angle=90, end angle=-38,radius=0.3];
    \node at (2, 2.25) {\scalebox{0.75}{$\theta_{i+1}$}};
    \node at (1, 2) [left] {\scalebox{0.75}{$z_{i+1}$}};
    \node at (2, 1.75) {\scalebox{0.75}{$\alpha_i+\beta_i$}};
    \end{tikzpicture}
    \caption{Illustration of relative coordinates.}\label{fig-rel.coordin.}
    \end{center}
\end{figure}
\\
As noted in~\cite{marshall2004formations}, the transformation from $q_i=\left[x_i,y_i,\theta_i\right]^\top$ into $\xi_i=\left[\rho_i,\alpha_i,\beta_i\right]^\top$ is not invertible due to the removal of any reference to a global coordinate frame. In addition, it gives rise to some constraints on the relative coordinates.
\\
Specifically, $\sum_{i=1}^n{z_{i+1}-z_i}\equiv0$ holds over time. Let $\xi=\left[\xi_1^\top,\xi_2^\top,\ldots,\xi_n^\top\right]^\top$ represent the full system state in relative coordinates. By choosing a coordinate frame attached to (say) unicycle 1 and oriented with its heading direction, the above identity imposes the following constraints:
\begin{align}
    &g_1(\xi) \coloneqq \rho_1\sin{\alpha_1}+\rho_2\sin{(\alpha_2+\pi-\beta_1)}+\ldots \nonumber \\
    &~ +\rho_n\sin{(\alpha_n+(n-1)\pi-\beta_1-\ldots-\beta_{n-1})} = 0, \label{eq-constr-1} \\
    &g_2(\xi) \coloneqq \rho_1\cos{\alpha_1}+\rho_2\cos{(\alpha_2+\pi-\beta_1)}+\ldots \nonumber \\
    &~ +\rho_n\cos{(\alpha_n+(n-1)\pi-\beta_1-\ldots-\beta_{n-1})} = 0. \label{eq-constr-2} 
\end{align}
Fig.~\ref{fig-geomet.constr.} illustrates how these constraints arise: In the coordinate frame attached to (say) unicycle 1 and oriented with its heading direction, the horizontal (vertical) projections of the distances $\rho_i$'s sum to zero. 
\begin{figure}[h]
    \centering
    \begin{tikzpicture}
    \draw[->, thick] (0,0) -- (1,0);
    \node at (0,0) [left] {\scalebox{0.75}{$1$}};
    \draw[dashed, thick] (0,0) -- (1.2,{1.2*sqrt(3)}) node[midway, fill=white, inner sep=1pt] {\scalebox{0.75}{$\rho_1$}};
    \draw[->, thick] (1.2,{1.2*sqrt(3)}) -- (1.2,{1.2*sqrt(3)+1});
    \node at (1.3,{1.1*sqrt(3)}) {\scalebox{0.75}{$2$}};
    \draw[dashed, thick] (1.2,{1.2*sqrt(3)}) -- (1.7,{1.7*sqrt(3)});
    \draw[dashed, thick] (1.2,{1.2*sqrt(3)}) -- (2.2,{1.2*sqrt(3)});
    \draw[dashed, thick] (1.2,{1.2*sqrt(3)}) -- ({1.2-sqrt(3))},{1+1.2*sqrt(3)}) node[midway, fill=white, inner sep=1pt] {\scalebox{0.75}{$\rho_2$}};
    \draw[->, thick] ({1.2-sqrt(3))},{1+1.2*sqrt(3)}) -- ({1.2-1.5*sqrt(3)},{.7+1.2*sqrt(3)});
    \node at ({1.4-sqrt(3))},{1+1.25*sqrt(3)}) {\scalebox{0.75}{$3$}};
    \draw[dashed, thick] ({1.2-sqrt(3)},{1+1.2*sqrt(3)}) -- ({.7-sqrt(3))},{1+.7*sqrt(3)});

    \draw[->, thick] (.5,0) arc[start angle=0, end angle=60,radius=0.5];
    \node at (.75,.25) {\scalebox{0.75}{$\alpha_1$}};
    \draw[->, thick] (1.7,{1.2*sqrt(3)}) arc[start angle=0, end angle=60,radius=0.5];
    \node at (1.95,{1.2*sqrt(3)+.25}) {\scalebox{0.75}{$\alpha_1$}};
    \draw[->, thick] (1.2,{1.2*sqrt(3)+.5}) arc[start angle=90, end angle=150,radius=0.5];
    \node at (.9,{1.2*sqrt(3)+.6}){\scalebox{0.75}{$\alpha_2$}};
    \draw[->, thick] (1.2,{1.2*sqrt(3)+.3}) arc[start angle=90, end angle=240,radius=0.3];
    \node at (.4,{1.2*sqrt(3)}){\scalebox{0.75}{$\alpha_1+\beta_1$}};
    \draw[->, thick] (1.55,{1.55*sqrt(3)}) arc[start angle=60, end angle=90,radius=0.7];
    \node at ({1.2*sqrt(3)},{1.2*sqrt(3)+1}) {\scalebox{0.75}{$\pi-\alpha_1-\beta_1$}};

    \end{tikzpicture}
    \caption{Illustration of coordinate constraints.}
    \label{fig-geomet.constr.}
\end{figure} 
\\
Moreover, $\sum_{i=1}^n{\dot{\beta}_i(t)}=0~\Rightarrow~\sum_{i=1}^n{\beta_i(t)}\equiv c,~\forall t\geq0$, where $c=-n\pi$ follows from the definition of  $\beta_i$, leading to the final coordinate constraint: 
\begin{equation}
    g_3(\xi) \coloneqq \sum_{i=1}^n{\beta_i+n\pi} = 0~\text{mod}~2\pi. \label{eq-constr-3} 
\end{equation}
These constraints are always satisfied regardless of the control design, a fact that will be elaborated in the next section. 

\subsection{Circular formation}\label{sec-equilibria}
In this part we develop a systematic analysis of circular formation. A necessary and sufficient condition for the existence of circular formation in cyclic pursuit is derived, which characterizes the target equilibrium formation associated with~\eqref{eq-relative-kinematic-1}. In addition, several properties satisfied in circular formation are discussed.
\begin{definition}
    Consider a group of unicycles without overlap, a circular formation arises if they travel along a common circle at the same angular speed.
\end{definition}
We first present the existence condition in global coordinates, which serves as a prerequisite for Theorem~\ref{thm-stable-criterion}, the main theorem of the paper. We then provide an equivalent formulation in relative coordinates, which provides insights into the control design.
\begin{proposition}\label{thm-exist.-condt.-global}
Consider a group of unicycles without overlap, a circular formation is realized if and only if the following conditions hold:
\begin{align}
    &\tilde{z}_j\times\tilde{\mathtt{r}}_j=0,~\forall j\in\{2,3,\ldots,n\}, \label{eq-product-zero} \\
    &\frac{\tilde{z}_2^\top\tilde{\mathtt{r}}_2}{\|\tilde{\mathtt{r}}_2\|^2} - \frac{\tilde{z}_j^\top\tilde{\mathtt{r}}_j}{\|\tilde{\mathtt{r}}_j\|^2} = 0,~\forall j\in\{3,4,\ldots,n\}, \label{eq-diff-zero}
\end{align}
where $\|\tilde{z}_j\|$ and $\|\tilde{\mathtt{r}}_j\|$ are nonzero for all $j\in\{2,3,\ldots,n\}$, and
\begin{align*}
    \tilde{z}_j &=z_1-z_j, ~z_*=\left[x_*,y_*\right]^\top, ~*\in\{1,j\}, \\
    \tilde{\mathtt{r}}_j &=\mathtt{r}_1-\mathtt{r}_j, ~\mathtt{r}_*=[-\sin{\theta_*},\cos{\theta_*}]^\top.
\end{align*}
The submanifold induced by \eqref{eq-product-zero}--\eqref{eq-diff-zero} has dimension $n+3$. 
\end{proposition}
The proof is placed in Appendix~\ref{app-exist.-condt.-global}.

Note that $\|\cdot\|$ denotes the $2$-norm, and $\times$ the cross product. Since $\|\tilde{z}_j\|$ and $\|\tilde{\mathtt{r}}_j\|$ are nonzero, $\tilde{z}_j\times\tilde{\mathtt{r}}_j=0\Leftrightarrow \tilde{z}_j=r_j\tilde{\mathtt{r}}_j$, where $r_j={\tilde{z}_j^\top\tilde{\mathtt{r}}_j}/{\|\tilde{\mathtt{r}}_j\|^2}\neq0$. This condition ensures that agents $1$ and $j\geq2$ travel along a circle, counterclockwise (resp. clockwise) for $r_j<0$ (resp. $r_j>0$). Condition~\eqref{eq-diff-zero} further ensures that $r\coloneqq r_2=r_j$ for all $j\geq3$, and thus all agents travel along a common circle. An illustration of~\eqref{eq-product-zero} is provided in Fig.~\ref{fig-product-zero}.
\begin{figure}
    \centering
    \begin{tikzpicture}

    \draw[->, thick] (-2,0) -- (-1.5,0) node[right] {\scalebox{0.75}{$x$}};
    \draw[->, thick] (-2,0) -- (-2,.5) node[above] {\scalebox{0.75}{$y$}};
    \node at (2,0)  [right] {\scalebox{0.75}{$z_\mathcal{O}$}};
    
    \draw[->, thick] (0,0) -- (0, -.75) node[below] {\scalebox{0.75}{$\mathtt{v}_1$}};
    \node at (0,0) [left] {\scalebox{0.75}{$z_1$}};
    \draw[->, thick] (.15,0) -- (.75, 0) node[below] {\scalebox{0.75}{$\mathtt{r}_1$}};
    
    \draw[->, thick] (2,2) -- (1.25,2) node[left] {\scalebox{0.75}{$\mathtt{v}_j$}};
    \node at (2, 2) [right] {\scalebox{0.75}{$z_j$}};
    \draw[->, thick] (2,2) -- (2,1.25) node[right] {\scalebox{0.75}{$\mathtt{r}_j$}};

    \draw[thick] (2-.2, 2) -- (2-.2, 2-.2);
    \draw[thick] (2, 2-.2) -- (2-.2, 2-.2);
    \draw[thick] (0, -.2) -- (.2, -.2);
    \draw[thick] (.2, 0) -- (.2, -.2);

    \draw[dashed, thick] (.85,.0) -- (2,0);
    \draw[dashed, thick] (2,0) -- (2,.65);
    \draw[dashed, thick] (2,.85) -- (2,1.2);

    \draw[->, thick] (2,2) -- (0,0) node[midway, fill=white, inner sep=1pt] {\scalebox{0.75}{$\tilde{z}_j$}};
    \draw[->, thick] (1.25,0) -- (2,.75) node[midway, fill=white, inner sep=1pt] {\scalebox{0.75}{$\tilde{\mathtt{r}}_j$}};
    \draw[dashed, thick] (2,2) arc[start angle=90, end angle=173, radius=2];

    \end{tikzpicture}
    \caption{Illustration of $\tilde{z}_j\times\tilde{\mathtt{r}}_j=0$.}
    \label{fig-product-zero}
\end{figure}

\begin{proposition}\label{thm-exist.-condt.}
    Consider a cyclic pursuit of unicycles without overlap, a circular formation is realized if and only if the following extra conditions hold (in addition to~\eqref{eq-constr-1}–\eqref{eq-constr-3}):
    \begin{align}
        2\alpha_i+\beta_i&=\pi~\text{mod}~2\pi,~\forall i\in\mathcal{I},\label{eq-sum-pi-mod-2pi} \\
        {2\sin{\alpha_i}}/{\rho_i} &= {1}/{r}\neq0,~\forall i\in\mathcal{I},\label{eq-same-ratio}
    \end{align}
    where $|r|\in(0,+\infty)$ is equal to the radius of the circle. The unicycles rotate at angular speed $v/r$, counterclockwise (resp. clockwise) for $r>0$ (resp. $r<0$).
\end{proposition}
The proof is placed in Appendix~\ref{app-exist.-condt.}.

\begin{remark}\label{rmk-dimension-dismatch}
    With all $\sin{\alpha_i}$'s being nonzero, we can reformulate~\eqref{eq-same-ratio} as ${\sin{\alpha_i}}/{\rho_i}-{\sin{\alpha_{i+1}}}/{\rho_{i+1}}= 0$.  Since the  constraints~\eqref{eq-constr-1}–\eqref{eq-constr-3} are always satisfied, there are three redundant coordinates (say) $\rho_{n}$, $\alpha_n$, and $\beta_n$. Thus, it suffices to consider $n-1$ equations from~\eqref{eq-sum-pi-mod-2pi} for $i\in\mathcal{I}\setminus\{n\}$ and $n-2$ equations from~\eqref{eq-same-ratio} for $i\in\mathcal{I}\setminus\{n-1,n\}$, yielding a total of $2n-3$ equations that defines a submanifold corresponding to circular formation, which has dimension $n$ in the relative-coordinate space (of dimension $3n-3$), rather than $n+3$ in the global-coordinate space. The three missing dimensions represent the translation and rotation of the formation in global coordinates.
\end{remark}

Condition~\eqref{eq-sum-pi-mod-2pi} (resp. \eqref{eq-same-ratio}) is a counterpart of condition~\eqref{eq-product-zero} (resp. \eqref{eq-diff-zero}) in relative coordinates. An illustration of~\eqref{eq-sum-pi-mod-2pi} is provided in Fig.~\ref{fig-sum-pi}. A corollary of~\eqref{eq-sum-pi-mod-2pi} is derived, whose proof is straightforward and thus omitted.
\begin{figure}
    \centering
    \begin{tikzpicture}

    \draw[->, thick] (-2,0) -- (-1.5,0) node[right] {\scalebox{0.75}{$x$}};
    \draw[->, thick] (-2,0) -- (-2,.5) node[above] {\scalebox{0.75}{$y$}};
    \node at (2,0)  [right] {\scalebox{0.75}{$z_\mathcal{O}$}};
    
    \draw[->, thick] (0,0) -- (0, -.75) node[below] {\scalebox{0.75}{$\mathtt{v}_{i+1}$}};
    \node at (0,0) [left] {\scalebox{0.75}{$z_{i+1}$}};
    \draw[->, thick] (0,0) -- (.75, 0) node[right] {\scalebox{0.75}{$\mathtt{r}_{i+1}$}};
    \draw[dashed, thick] (1.5,0) -- (2,0);
    \draw[->, thick] (0, -.4) arc[start angle=-90, end angle=45, radius=0.4];
    \node at (.8, -.4) {\scalebox{0.75}{$\alpha_i+\beta_i$}};
    
    \draw[->, thick] (2,2) -- (1.25,2) node[left] {\scalebox{0.75}{$\mathtt{v}_i$}};
    \node at (2, 2) [right] {\scalebox{0.75}{$z_i$}};
    \draw[dashed, thick] (2,.8) -- (2,0);
    \draw[->, thick] (2,2) -- (2,1.25) node[below] {\scalebox{0.75}{$\mathtt{r}_i$}};
    \draw[->, thick] (2-.4,2) arc[start angle=180, end angle=225, radius=0.4];
    \node at (1.45, 2-.25) {\scalebox{0.75}{$\alpha_i$}};

    \draw[thick] (2-.2, 2) -- (2-.2, 2-.2);
    \draw[thick] (2, 2-.2) -- (2-.2, 2-.2);
    \draw[thick] (0, -.2) -- (.2, -.2);
    \draw[thick] (.2, 0) -- (.2, -.2);

    \draw[thick] (2,2) -- (0,0) node[midway, fill=white, inner sep=1pt] {\scalebox{0.75}{$\rho_i$}};
    \draw[dashed, thick] (2,2) arc[start angle=90, end angle=180, radius=2];

    \end{tikzpicture}
    \caption{Illustration of $2\alpha_i+\beta_i= \pi$, corresponding to counterclockwise rotation.}
    \label{fig-sum-pi}
\end{figure}

\begin{corollary}\label{thm-exist.-condt.-equiv}
Condition~\eqref{eq-sum-pi-mod-2pi} is equivalent to the following equation, which must hold for all $i\in\mathcal{I}$:
    \begin{subequations}\label{eq-sum-pmpi-2}
        \begin{align}
            \sin{\alpha_i}-\sin{(\alpha_i+\beta_i)} &= 0, \\
            \cos{\alpha_i}+\cos{(\alpha_i+\beta_i)} &= 0.
        \end{align} 
    \end{subequations}
\end{corollary}
\begin{remark}
    Without loss of generality, in the rest of the paper we normalize all angle variables, including $\alpha_i+\beta_i$, to the range $[-\pi,\pi)$ whenever necessary.
    This normalization simplifies the analysis and is justified by the fact that both~\eqref{eq-relative-kinematic-1} and the mechanism of circular formation depend on their trigonometric function values, rather than the angles themselves. 
\end{remark}
In what follows, we present a set of definitions and propositions that are fundamental to our equilibrium and stability analysis. 
\begin{definition}\label{def-arrangements}
    The arrangement of a circular formation  refers to the relative ordering of the unicycles along the circle in the rotation direction.
\end{definition}
\begin{proposition}\label{thm-arrangements}
    For $n$ unicycles in a circular formation, there are $(n-1)!$ possible distinct arrangements.
\end{proposition}
\proof
    There are $(n-1)!$ distinct permutations when arranging $n$ distinct points on a circle.
\endproof
\begin{figure}[h]
    \centering
    {
        \begin{tikzpicture}[scale = 0.75]
            \draw[->, thick] (0,1.5) -- (-.75,1.5) node[left] {\scalebox{0.75}{$\mathtt{v}_1$}};
            \draw[->, thick] (-1.5,0) -- (-1.5,-.75)  node[below] {\scalebox{0.75}{$\mathtt{v}_2$}};
            \draw[->, thick] ({1.5/sqrt(2)},{-1.5/sqrt(2)}) -- ({1.5/sqrt(2)+.75/sqrt(2)}, {-1.5/sqrt(2)+.75/sqrt(2)}) node[right] {\scalebox{0.75}{$\mathtt{v}_3$}};
            
            \draw[->, thick, >=latex, line width=1pt] (0,1.5) -- (-1.5,0);
            \draw[->, thick, >=latex, line width=1pt] (-1.5,0) -- ({1.5/sqrt(2)},{-1.5/sqrt(2)});
            \draw[->, thick, >=latex, line width=1pt] ({1.5/sqrt(2)},{-1.5/sqrt(2)}) -- (0,1.5);
        
            \draw[dashed, thick] (0,0) -- (0,1.5);
            \draw[dashed, thick] (0,0) -- (-1.5,0);
            \draw[dashed, thick] (0,0) -- ({1.5/sqrt(2)},{-1.5/sqrt(2)});
            
            \draw[->, thick, >=latex, line width=1pt] (0:1.5) arc[start angle=0, end angle=45, radius=1.5];
            
            \draw[dashed, thick] (45:1.5) arc[start angle=45, end angle=360, radius=1.5];
        \end{tikzpicture}
    }
    \hspace{0.1\linewidth}
    {
        \begin{tikzpicture}[scale = 0.75]
            \draw[->, thick] (0,1.5) -- (-.75,1.5) node[left] {\scalebox{0.75}{$\mathtt{v}_1$}};
            \draw[->, thick] (1.5,0) -- (1.5,.75)  node[above] {\scalebox{0.75}{$\mathtt{v}_2$}};
            \draw[->, thick] ({-1.5/sqrt(2)},{-1.5/sqrt(2)}) -- ({-1.5/sqrt(2)+.75/sqrt(2)}, {-1.5/sqrt(2)-.75/sqrt(2)}) node[left] {\scalebox{0.75}{$\mathtt{v}_3$}};
            
            \draw[<-, thick, >=latex, line width=1pt] (0,1.5) -- ({-1.5/sqrt(2)},{-1.5/sqrt(2)});
            \draw[<-, thick, >=latex, line width=1pt] ({-1.5/sqrt(2)},{-1.5/sqrt(2)}) -- (1.5,0);
            \draw[<-, thick, >=latex, line width=1pt] (1.5,0) -- (0,1.5);
        
            \draw[dashed, thick] (0,0) -- (0,1.5);
            \draw[dashed, thick] (0,0) -- (1.5,0);
            \draw[dashed, thick] (0,0) -- ({-1.5/sqrt(2)},{-1.5/sqrt(2)});

            \draw[dashed, thick] (0:1.5) arc[start angle=0, end angle=135, radius=1.5];

            \draw[dashed, thick] (180:1.5) arc[start angle=180, end angle=360, radius=1.5];

            \draw[->, thick, >=latex, line width=1pt] (135:1.5) arc[start angle=135, end angle=180, radius=1.5];
        \end{tikzpicture}
    }
    \caption{Two arrangements of counterclockwise formations with three vehicles. The second formation is obtained by reversing the rotation direction in the first one and viewing from inside the paper outward.}
    \label{fig-geomet.-config.-n=3}
\end{figure} 

Fig.~\ref{fig-geomet.-config.-n=3} gives an illustration of arrangements using $n=3$ as an example. The next proposition claims that $\sum_{i=1}^n{\alpha_i}$ is a nonzero multiple of $\pi$ when $n$ unicycles form a circular formation. Note that $\mathbb{Z}_+$ denotes the set of positive integers. If $\alpha_i$'s are not restricted to $[-\pi,\pi)$, the generalized result $\sum_{i=1}^n{\alpha_i}=0~\text{mod}~\pi$ holds.
\begin{proposition}\label{thm-alpha-sum}
    For any circular formation, the sum of bearings satisfies 
    \begin{equation}
        \sum_{i=1}^n{\alpha_i} = \pm\mathtt{p}\pi, \quad \text{where}\quad \mathtt{p} \in \left[1,(n-1)\right] \subset \mathbb{Z}_+. \label{eq-sum-multiple} 
    \end{equation}
    The sum is positive (resp. negative) for counterclockwise (resp. clockwise) rotation. 
\end{proposition}
The proof is placed in Appendix~\ref{app-alpha-sum}.

\begin{definition}\label{def-regular-motion}
    A circular formation is said to be regular if $\mathtt{p}=1$; otherwise, it is said to be irregular.
\end{definition}
A formation with $\mathtt{p}=1$ is termed regular, in the sense that, following the rotation direction, the unicycles are arranged along the circle in cyclic order. 
\begin{proposition}\label{thm-cot-sum}
    For any counterclockwise formation with $n\geq3$, it holds that
    \begin{equation}
        \sum_{i=1}^n{\cot{\alpha_i}}\in\left\{\begin{aligned}
            &(0,\infty),~\text{if}~\mathtt{p} = 1; \\
    &(-\infty,0),~\text{if}~\mathtt{p} = n-1; \\
    &\mathbb{R},~\text{otherwise}.
        \end{aligned}\right.
    \end{equation}
\end{proposition}
The proof is placed in Appendix~\ref{app-cot-sum}.

Reversed results for clockwise formation can be obtained due to $\cot{(-x)} = -\cot{x}$. When $n=2$, it is trivial that $\sum_{i=1}^n{\cot{\alpha_i}}=0$.

\section{Control design and stability analysis}\label{sec-controller}
In this section we propose a local control law that enables the emergence of circular formation. We then explore the spectrum of the linearized closed-loop system matrix, revealing key features that facilitate stability analysis and model reduction.
\subsection{Emergence-oriented control}
The following control law with a gain parameter $k\neq0$ is proposed to achieve circular formation: 
\begin{equation}\label{eq-omega-1}
    \omega_i = \frac{2v}{\rho_i}\sin{\alpha_i}+\frac{k}{\rho_i}\left(\cos{\alpha_i}+\cos{(\alpha_i+\beta_i)}\right).
\end{equation}
Substituting~\eqref{eq-omega-1} into~\eqref{eq-relative-kinematic-1} gives rise to a cyclically interconnected system of identical nonlinear subsystems:
\begin{subequations}\label{eq-relative-kinematic-2}\begin{align}
    \dot{\rho_i} &= -v\left(\cos{\alpha_i}+\cos{(\alpha_i+\beta_i)}\right), \label{eq-relative-kinematic-2-rho} \\
    \dot{\alpha_i} &= \frac{v}{\rho_i}\left( \sin{(\alpha_i+\beta_i)}-\sin{\alpha_i}\right)  -\frac{k}{\rho_i}\left(\cos{\alpha_i}+\cos{(\alpha_i+\beta_i)}\right), \label{eq-relative-kinematic-2-alpha} \\
    \dot{\beta}_i &= \frac{2v}{\rho_i}\sin{\alpha_i} -\frac{2v}{\rho_{i+1}}\sin{\alpha_{i+1}} + \frac{k}{\rho_i}\left(\cos{\alpha_i}+\cos{(\alpha_i+\beta_i)}\right) \nonumber \\ 
    &\quad  -\frac{k}{\rho_{i+1}}\left(\cos{\alpha_{i+1}}+\cos{(\alpha_{i+1}+\beta_{i+1})}\right). \label{eq-relative-kinematic-2-beta} 
\end{align}\end{subequations}
View these subsystems as $\dot{\xi}_i=f(\xi_i,\xi_{i+1})$ and the full system as 
\begin{equation}\label{eq-complete-syst}
    \dot{\xi}=\hat{f}(\xi).
\end{equation}
Along the lines of~\cite{marshall2004formations}, we show that the system is constrained to evolve on a $\hat{f}$-invariant submanifold of $\mathbb{R}^{3n}$, which explains why the three constraints are always satisfied. This invariance is key to understanding how the spectrum of the linearized system matrix relates to the stability of the target equilibrium formation. 
\\
Let $g(\xi)=\left[g_1(\xi),g_2(\xi),g_3(\xi)\right]^\top\in\mathbb{R}^3$, then
\begin{equation}
    \label{eq-manifold}
    \mathcal{M} \coloneqq \left\{\xi\in\mathbb{R}^{3n}|g(\xi)=0\right\}\subset\mathbb{R}^{3n}
\end{equation}
defines a submanifold of dimension $3n-3$.
\begin{lemma}\label{thm-f-iinvariance}
    The submanifold $\mathcal{M}$ is invariant under $\hat{f}$.
\end{lemma}
The proof is placed in Appendix~\ref{app-f-iinvariance}.

\begin{theorem}\label{thm-equilibria}
    The invariant set (which in this paper is taken to be the set of equilibria) of the system $\dot{\xi}=\hat{f}(\xi)$ is given by 
    \begin{equation}
    \begin{aligned}
    \Xi \coloneqq \big\{
    &\bar{\xi}\in\mathcal{M} : 2\bar{\alpha}_i+\bar{\beta}_i=\pi~\text{mod}~2\pi \\
    &\text{and}~\exists\bar{s}~\text{s.t.}~{\sin{\bar{\alpha}}_i}/{\bar{\rho}_i}=\bar{s},~\forall i\in\mathcal{I}
    \big\}. \label{eq-set-equilibria}
    \end{aligned}\end{equation}
    The unicycles travel at angular speed $\bar{\omega}\coloneqq2v\bar{s}$. A circular formation is realized whenever $\bar{s} \neq 0$: counterclockwise for $\bar{s}>0$ and clockwise for $\bar{s}<0$. Otherwise, a collinear formation is realized; they travel along a line at the same forward speed.
\end{theorem}
\proof For $\dot{\rho}_i=0$ and $\dot{\alpha}_i=0$, \eqref{eq-relative-kinematic-2-rho}--\eqref{eq-relative-kinematic-2-alpha} gives $\cos{\bar{\alpha}_i}+\cos{(\bar{\alpha}_i+\bar{\beta}_i)}=0$ and $\sin{\bar{\alpha}_i}-\sin{(\bar{\alpha}_i+\bar{\beta}_i)}=0$, which are equivalent to the condition~\eqref{eq-sum-pi-mod-2pi} by Corollary~\ref{thm-exist.-condt.-equiv}. For $\dot{\beta}_i=0$, \eqref{eq-relative-kinematic-2-beta} yields $\sin{\bar{\alpha}_i}/\bar{\rho}_i = \sin{\bar{\alpha}_{i+1}}/\bar{\rho}_{i+1}$. Since~\eqref{eq-relative-kinematic-2-beta} is cyclically interconnected, $\sin{\bar{\alpha}_i}/\bar{\rho}_i=\bar{s}$ for all $i$, and $\omega_i=\bar{\omega}=2v\bar{s}$. Finally, $\bar{\xi}\in\mathcal{M}$ as $\mathcal{M}$ is $\hat{f}$-invariant. 
\\
When $\bar{\omega}=2v\bar{s}\neq0$, each equilibrium corresponds to a particular circular formation of diameter $1/{|\bar{s}|}$, with the agent spacing and rotation direction determined accordingly; see~Proposition~\ref{thm-exist.-condt.}.  Otherwise, it follows that $\bar{\alpha}_i=0~\text{mod}~\pi$ and $\bar{\beta}_i=-\pi~\text{mod}~2\pi$ for all $i$, implying a collinear formation.
\endproof

\subsection{Stability analysis of invariant sets}
Write
$\bar{\xi}=\left[\bar{\xi}_1^{\;\top},\bar{\xi}_2^{\;\top},\ldots,\bar{\xi}_n^{\;\top}\right]^\top\in\Xi$ with $\bar{\xi}_i=\left[\bar{\rho}_i,\bar{\alpha}_i,\bar{\beta}_i\right]^\top$, and $\tilde{\xi}_i=\xi_i-\bar{\xi}_i$, then linearizing $\tilde{\xi}_i$ at any $\bar{\xi}$ gives rise to $n$ identical subsystems of the form $\dot{\tilde {\xi}}_i=A_i\tilde{\xi}_i+B_{i+1}\tilde{\xi}_{i+1}$. Denote $\bar{c}_i \coloneqq {\cos{\bar{\alpha}_i}}/{\bar{\rho}_i}$, then we have
\begin{equation*}\begin{aligned}
    A_i &= \begin{bmatrix}
        0 & 2v\bar{s}\bar{\rho}_i & v\bar{s}\bar{\rho}_i \\
        0 & 2k\bar{s}-2v\bar{c}_i & k\bar{s}-v\bar{c}_i \\
        -2v\bar{s}/\bar{\rho}_i & 2v\bar{c}_i-2k\bar{s} & -k\bar{s}
    \end{bmatrix}, \\
    B_{i+1} &= \begin{bmatrix}
        0 & 0 & 0 \\
        0 & 0 & 0 \\
        2v\bar{s}/\bar{\rho}_{i+1} & 2k\bar{s}-2v\bar{c}_{i+1} & k\bar{s}
    \end{bmatrix}.
\end{aligned}\end{equation*}
Denote $\tilde{\xi}=\xi-\bar{\xi}$, then the full linearized system has the form $\dot{\tilde{\xi}}=\hat{A}\tilde{\xi}$, where
\begin{align}\label{eq-A-hat}
    \hat{A} = \begin{bmatrix}
        A_1 & B_2 \\
        & A_2 & B_3 \\
        & & \ddots & \ddots \\
        B_1 & & &  A_n
    \end{bmatrix}.
\end{align}
Collinear formation is the limiting case of circular formation with infinite radius. While the formation of collinear formations forms a feasible invariant set, we show that it is not asymptotically stable. Hence, its stability warrants no further analysis, and we henceforth focus on the formation of circular formations.
\begin{theorem}\label{thm-unstable}
    The invariant set $\mathcal{L}\coloneqq\left\{\bar{\xi}\in\Xi : \bar{s}=0\right\}$ is not asymptotically stable, i.e., the formation of collinear formations is not asymptotically stable.
\end{theorem}
\proof For all $\bar{\xi}\in\mathcal{L}$, since $\bar{s}=0$ and all $\bar{c}_i$'s are nonzero, it follows that
\begin{equation*}
    \mathcal{S}_\mathsf{b}(\bar{\xi})
    \coloneqq \mspan{\left\{\left[0,\frac{1}{\bar{c}_1},-\frac{2}{\bar{c}_1},\ldots,0,\frac{1}{\bar{c}_n},-\frac{2}{\bar{c}_n}\right]^\top\right\}}\subset\ker{\hat{A}}.
\end{equation*}
$\mathcal{S}_\mathsf{b}(\bar{\xi})$ characterizes infinitesimal motions that deform a line into an arc, which appears nearly straight at small scales. For any nonzero $\tilde{\xi} \in \mathcal{S}_\mathsf{b}(\bar{\xi})$, the perturbed state $\bar{\xi} + \tilde{\xi}\notin\mathcal{L}$; therefore, $\mathcal{L}$ is not asymptotically stable.  
\endproof 

Recall that the submanifold $\mathcal{M}$ is $\hat{f}$-invariant, so there exists an induced linear transformation in the quotient space $\mathbb{R}^{3n}/T_{\bar{\xi}}\mathcal{M}$ whose eigenvalues do not influence the stability of $\bar{\xi}\in\mathcal{M}$, where $T_{\bar{\xi}}\mathcal{M}$ denotes the tangent space of $\mathcal{M}$ at $\bar{\xi}\in\mathcal{M}$. Thus, there exists a change of coordinates for $\mathbb{R}^{3n}$ that transforms $\hat{A}$ into block triangular form
\begin{equation}\label{eq-A-upper-form}
    \begin{bmatrix}
        \hat{A}_{T_{\bar{\xi}}\mathcal{M}} & * \\ 0_{3\times(3n-3)} & \hat{A}_{T_{\bar{\xi}}\mathcal{M}}^\star
    \end{bmatrix}.
\end{equation}
\begin{lemma}\label{thm-A-hat-pmjomega}
    In the quotient space $\mathbb{R}^{3n}/T_{\bar{\xi}}\mathcal{M}$, the induced linear transformation $\hat{A}_{T_{\bar{\xi}}\mathcal{M}}^\star$ has solely imaginary eigenvalues $\lambda_{1,2,3} = \{0,\pm j\bar{\omega}\}$.
\end{lemma}
The proof is placed in Appendix~\ref{app-A-hat-pmjomega}.

Therefore, when determining the stability of $\bar{\xi}\in\mathcal{M}$ we can disregard these three solely imaginary eigenvalues and analyze stability based on the remaining eigenvalues of $\hat{A}$. \\
On the other hand, a direct examination shows that the original $\hat{A}$ has $n+1$ distinct zero eigenvalues. They stem from the system’s degrees of freedom on the targeted equilibrium formation
\begin{equation}
    \mathcal{C} \coloneqq \left\{\bar{\xi}\in\Xi : \bar{s}\neq0\right\} = \Xi\setminus\mathcal{L},
\end{equation}
which allows circular formations with arbitrary radius and agent spacing. As a consequence, individual equilibria are not asymptotically stable, yet the stability of any neighborhood formed by a continuum of equilibria is unaffected.
\\
The invariant set $\mathcal{C}$ is disconnected, and naturally splits into two disjoint subsets:
\begin{subequations}\label{eq-C-two-parts}\begin{align}
    \mathcal{C}^-&\coloneqq\left\{\bar{\xi}\in\mathcal{C} : k\bar{s}<0\right\}, \label{eq-C-minus} \\
    \mathcal{C}^+&\coloneqq\left\{\bar{\xi}\in\mathcal{C} : k\bar{s}>0\right\} \label{eq-C-plus}.
\end{align}\end{subequations}
Notice that $\mathcal{C}^-$ and $\mathcal{C}^+$ are disconnected for $n>2$. By Propositions~\ref{thm-arrangements}--\ref{thm-alpha-sum}, $\mathcal{C}^-$ decomposes into a finite number of disjoint components (i.e., maximal connected subsets) of the form
\begin{equation}\label{eq-C-of-pi}
    \mathcal{C}^-_{\mathtt{p}\pi}\coloneqq\left\{\bar{\xi}\in\mathcal{C}^- : \sum_{i=1}^n{|\bar{\alpha}_i|=\mathtt{p}\pi}\right\},
\end{equation}
where $\mathtt{p}\in[1,(n-1)]\subset \mathbb{Z}_+$ and we write $|\bar{\alpha}_i|$ instead of $\bar{\alpha}_i$ to cover both counterclockwise and clockwise rotations; an analogous decomposition holds for $\mathcal{C}^+$. Since the stability may differs across these subsets (as discussed in the next subsection), the stability analysis is conducted with respect to any maximal connected subset of interest.
\\ 
The preceding discussion is then formalized and further developed in the following theorem.
\begin{theorem}\label{thm-stable-criterion}
Let $\mathcal{C} \coloneqq \left\{\bar{\xi}\in\Xi : \bar{s}\neq0\right\}$ and $\mathcal{C}^*\subset\mathcal{C}$ be a maximal connected subset.
(i) For all $\bar{\xi}\in\mathcal{C}$, the matrix $\hat{A}$ has a pair of purely imaginary eigenvalues $\pm j2v\bar{s}$ and $n+1$ zero eigenvalues. 
(ii) If, for all $\bar{\xi}\in\mathcal{C}^*$, the remaining $2n-3$ eigenvalues have negative real parts, then $\mathcal{C}^*$ is locally asymptotically stable, i.e., there exists sufficiently small $\delta>0$ such that
\begin{equation}\begin{aligned}
    \xi(0)\in \mathcal{U}_\delta(\mathcal{C}^*) &\Rightarrow
    \lim_{t\rightarrow\infty}{\dist{\left(\xi(t),\mathcal{C}^*\right)}} = 0, \quad\text{where} \label{eq-set-stable} \\
    \mathcal{U}_\delta(\mathcal{C}^*)
    &\coloneqq\left\{\xi\in\mathcal{M}: \dist{\left(\xi,\mathcal{C}^*\right)}<\delta\right\}, \\
    \dist{\left(\xi,\mathcal{C}^*\right)}
    &\coloneqq \inf_{\bar{\xi}\in\mathcal{C}^*}{\left\|\xi-\bar{\xi}\right\|}.
\end{aligned}\end{equation}
\end{theorem}
The proof is placed in Appendix~\ref{app-set-stable}.

\subsection{Model reduction}\label{sec-reduction}
Whether the remaining $2n-3$ eigenvalues have negative real parts can be analyzed via the Routh–Hurwitz criterion applied to the characteristic polynomial. Prior to that,
we transform $\hat{A}$ into another block triangular form
\begin{equation}\label{eq-A-lower-form}
\begin{bmatrix}
    \hat{A}_{R} & 0_{2n\times n} \\
    * & 0_{n\times n} 
\end{bmatrix}.
\end{equation}
Thus, $\hat{A}_R$ is a $2n\times 2n$ matrix with solely imaginary eigenvalues $\{0,\pm j2v\bar{s}\}$. Then the problem reduces to computing the characteristic polynomial of $\hat{A}_R$, whose structure is detailed in the next lemma. 
\begin{lemma}\label{thm-A-reduced}
    For any $\bar{\xi}\in\mathcal{C}$, the matrix $\hat{A}$ can be transformed into the form~\eqref{eq-A-lower-form} via a similarity transformation, where $\hat{A}_R$ is given by
    \begin{equation}\label{eq-A-reduced}
        \hat{A}_R = \begin{bmatrix}
        \bar{A}_1 & \bar{B} \\
        & \bar{A}_2 & \bar{B} \\
        & & \ddots & \ddots \\
        \bar{B} & & & \bar{A}_n
    \end{bmatrix},
    \end{equation}
    where 
    $\bar{A}_i = \bar{s}\begin{bmatrix}
        0 & v\cot^2{\bar{\alpha}_i}-k\cot{\bar{\alpha}_i}+v \\
        -2v & k-2v\cot{\bar{\alpha}_i}
    \end{bmatrix}$ and $\bar{B} = \bar{s}\begin{bmatrix}
        0 & 0 \\  2v & k
    \end{bmatrix}$.
\end{lemma}
The proof is placed in Appendix~\ref{app-A-reduced}.

\begin{figure}
    \centering
    \includegraphics[width=.495\linewidth]{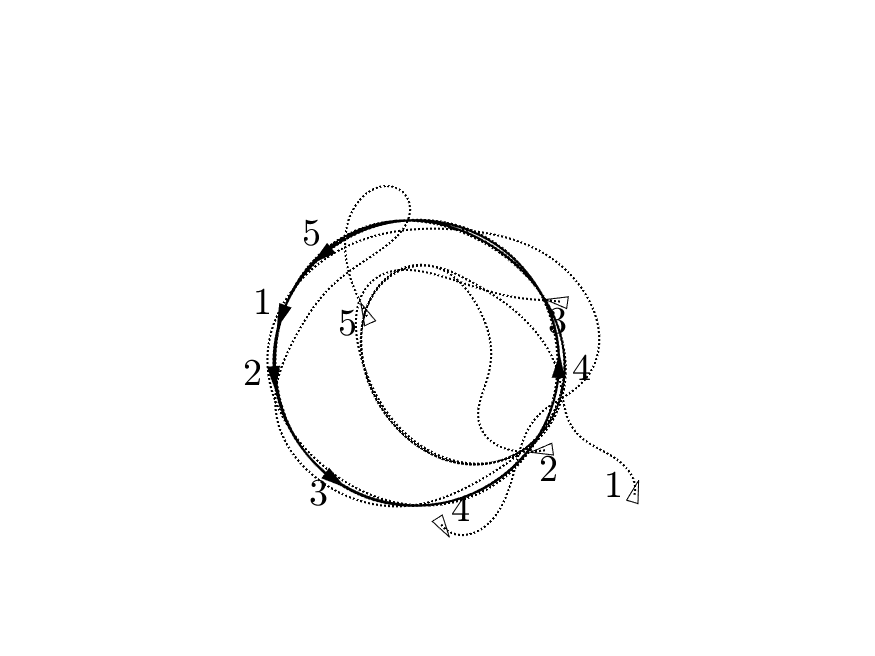}
    \includegraphics[width=.495\linewidth]{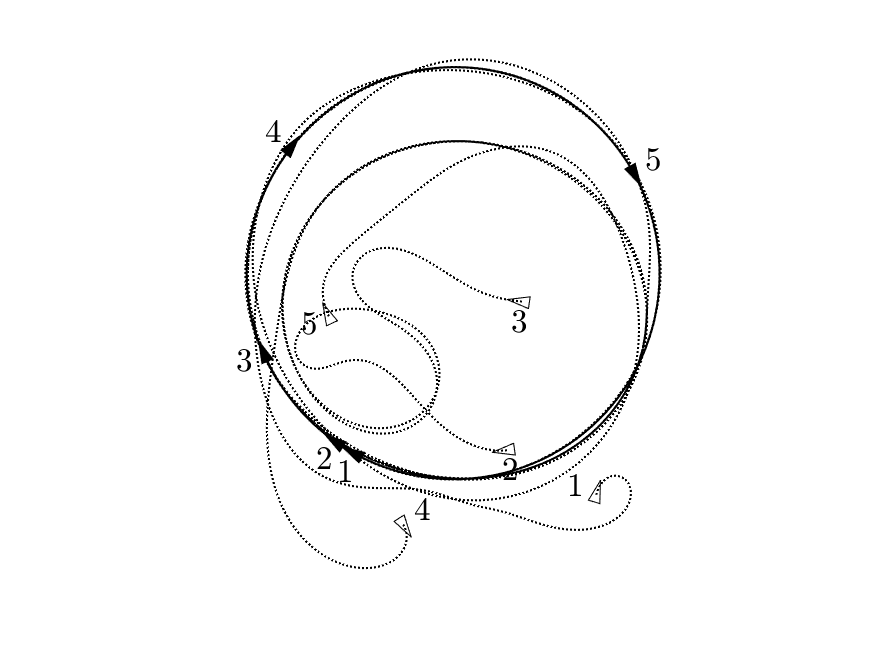}
    \caption{$n=5$ unicycles controlled by~\eqref{eq-omega-1} with $v=1$ and $k=\pm5$ converge to a circular formation. Under identical initial conditions, $k<0$ results in counterclockwise rotation (left), whereas $k>0$ results in clockwise rotation (right). Both formations are regular.}
    \label{fig-5vehicles}
\end{figure}

Although it yields a more pleasing structure, explicitly deriving the characteristic polynomial remains nontrivial. Instead, we derive an implicit expression. Before doing so, we discuss how the trace of $\hat{A}_R$ is influenced by $k$ and the different arrangements encoded in the bearing configurations. 
\\
We have $\tr{\hat{A}_R} =nk\bar{s}-2v\bar{s}\sum_{i=1}^n{\cot{\bar{\alpha}_i}}$. Since the trace of a matrix is equal to the sum of its eigenvalues, a negative trace is desirable; otherwise, either $\hat{A}_R$ is unstable (with eigenvalues in the open right-half plane), or it has all eigenvalues on the imaginary axis, rendering the linearization inconclusive.
\\
Since we wrap $\alpha_i$'s into the range $[-\pi,\pi)$, it holds that $\bar{s}\cdot{\cot{\bar{\alpha}_i}}=|\bar{s}|\cot{|\bar{\alpha}_i|}$. Thus, we can restrict our focus to counterclockwise rotation, i.e., $|\bar{\alpha}_i|\in(0,\pi)$. Rewrite the trace as 
\begin{equation*}
    \tr{\hat{A}_R} = |\bar{s}|\left(nk\cdot\sign{(\bar{s})}-2v\sum_{i=1}^n{\cot{|\bar{\alpha}_i|}}\right).
\end{equation*}
We can exclude $n$ and $2v$ from the discussion below since they are fixed. Then it follow from Proposition~\ref{thm-alpha-sum} that $\sum_{i=1}^n{|\bar{\alpha}_i|} = \mathtt{p}\pi$, with three possible cases depending on the value of $\mathtt{p}\in[1,n-1]\subset\mathbb{Z}_+$. \\
(i) For $\mathtt{p}=1$, i.e., regular arrangement, $-\sum_{i=1}^n{\cot{|\bar{\alpha}_i}|}$ is always negative. 
Thus, for any specific bearing configuration $\bar{\alpha}_i$'s, $k\cdot\sign{(\bar{s})}$ can be negative or positive, as long as the resulting trace is negative. \\
(ii) For $\mathtt{p}\in(1,n-1)$, the range of $-\sum_{i=1}^n{\cot{|\bar{\alpha}_i}|}$ is the whole $\mathbb{R}$. As a result, for some $\bar{\alpha}_i$'s (with $-\sum_{i=1}^n{\cot{|\bar{\alpha}_i}|}$ negative), $k\cdot\sign{(\bar{s})}$ can be negative or positive, while for others, it has to be negative. \\
(iii) For $\mathtt{p}=n-1$, which is the worst, $-\sum_{i=1}^n{\cot{|\bar{\alpha}_i}|}$ is always positive. Thus, in order that the trace is negative, $k\cdot\sign{(\bar{s})}$ has to be negative. \\
Given the variation in the feasible range of $k\cdot\sign{(\bar{s})}$ (so that the trace is negative) across different arrangements, it is more tractable to study their stability under a fixed $k$ rather than to find a full stabilizing range for each arrangement, \emph{if possible}. \\
Likewise, we focus solely on $\mathcal{C}^-$, i.e., $k\cdot\sign{(\bar{s})}<0$. The reason is twofold. First, it provides the most conservative setting in which (the stability of) all arrangements \emph{could} be analyzed. Second, numerical simulations suggest that counterclockwise (resp. clockwise) formation has a large region of attraction when $k<0$ (resp. $k>0$). An illustrative example is shown in~Fig.~\ref{fig-5vehicles}. 
\\
With the above lengthy yet necessary discussion in place, we are now ready to present the following theorem. Without loss of generality, let $v=1$. 
\begin{theorem}\label{thm-eigen.-poly.}
    Let $k=\pm2v$ and $v=1$. For any $\bar{\xi}\in\mathcal{C}^-$, an implicit expression for the characteristic polynomial of ${\hat{A}_R}/|\bar{s}|$ is
    \begin{align}
        P_{\bar{\xi}}(\lambda) &\coloneqq \prod_{i=1}^n{\left(\lambda^2+2z_i\lambda+2z_i^2\right)} -\prod_{i=1}^{n}{\left(2z_i^2-2\lambda\right)} \label{eq-implicit-poly.} \\
        &=\left(\lambda^3+4\lambda\right)\left(\lambda^{2n-3}+a_{2n-4}\lambda^{2n-4}+\ldots+a_0\right), \nonumber \\
        z_i &\coloneqq 1+\sign{(\bar{s})}\cot{\bar{\alpha}_i}. \label{eq-zi}
    \end{align}
\end{theorem}
The proof is placed in Appendix~\ref{app-eigen.-poly.}.

\begin{remark}\label{rmk-cot|alpha|}
    Since $\alpha_i$'s are normalized to $[-\pi,\pi)$, we can rewrite \eqref{eq-zi} as $z_i=1+\cot{|\bar{\alpha}_i|}$ with $|\bar{\alpha}_i|\in(0,\pi)$. Hence, we can focus solely on counterclockwise rotation.
\end{remark}

\section{Case studies of local stability}\label{sec-main} In this section we conclude a set of stability results for the case $n \leq 3$, which provides insights into the stability property for the case $n>3$. A special case is also studied, where $n>3$ unicycles are equally spaced in cyclic order. 
\subsection{Case Study: \texorpdfstring{$n=2$}{n=2}}
The favorable structure of $\hat{A}$ when $n=2$ allows us to determine the stability directly from its trace, as shown in the following theorem. 
\begin{theorem}\label{thm-stable-n=2}
    For $n=2$, the invariant set $\mathcal{C}^-$
    is locally asymptotically stable, whereas $\mathcal{C}^+$ is unstable, i.e., the formation of counterclockwise formations is locally asymptotically stable (resp. unstable) when $k<0$ (resp. $k>0$); the opposite holds for the formation of clockwise formations.
\end{theorem}
\proof We have $\tr{\hat{A}} =k\left(\bar{s}_1+\bar{s}_2\right)-2v\left(\bar{c}_1+\bar{c}_2\right) = 2k\bar{s}$. For every $\bar{\xi}\in\mathcal{C}$, $\hat{A}$ has five imaginary-axis eigenvalues whose sum is zero, so the last one is equal to the trace: $\lambda_6=2k\bar{s}\neq0$. The imaginary-axis eigenvalues can be disregarded. $\lambda_6>0$ for every $\bar{\xi}\in\mathcal{C}^+$, which follows that $\mathcal{C}^+$ is unstable. Conversely, $\lambda_6<0$ for every $\bar{\xi}\in\mathcal{C}^-$. Since $\mathcal{C}^-$ is a maximal connected subset (which is true only when $n=2$), it follows from Theorem~\ref{thm-stable-criterion} that $\mathcal{C}^-$ is locally asymptotically stable.
\endproof

\subsection{Case Study: \texorpdfstring{$n=3$}{n=3}}
In this part, the formation of regular formations is studied for $n=3$, and the following theorem is established. 
\begin{theorem}\label{thm-stable-n=3}
    For $n=3$ and $k=\pm2v$, the invariant set $\mathcal{C}_\pi^-$ is locally asymptotically stable, i.e.,
    the formation of regular counterclockwise (resp. clockwise) circular formations is locally asymptotically stable when $k=-2v$ (resp. $k=2v$).
\end{theorem}
The proof is placed in Appendix~\ref{app-stable-n=3}.

\begin{figure}
    \centering
    \includegraphics[width=.5\linewidth]{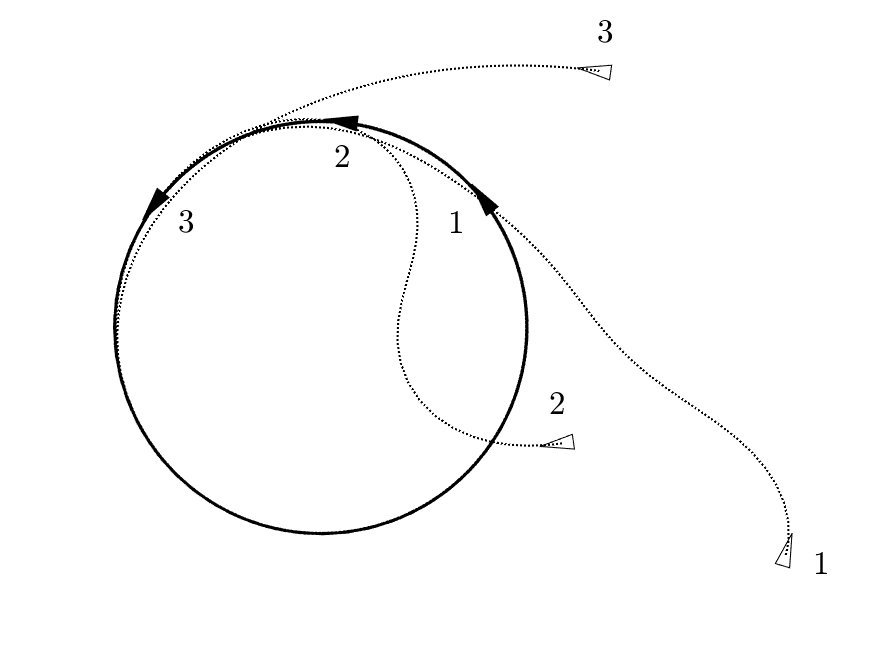}
    \caption{Three unicycles controlled by~\eqref{eq-omega-1} with $v=1$ and $k=-2$ converge to a regular counterclockwise formation.}
    \label{fig-3vehicles}
\end{figure}
\begin{figure}
    \centering
    \includegraphics[width=.5\linewidth]{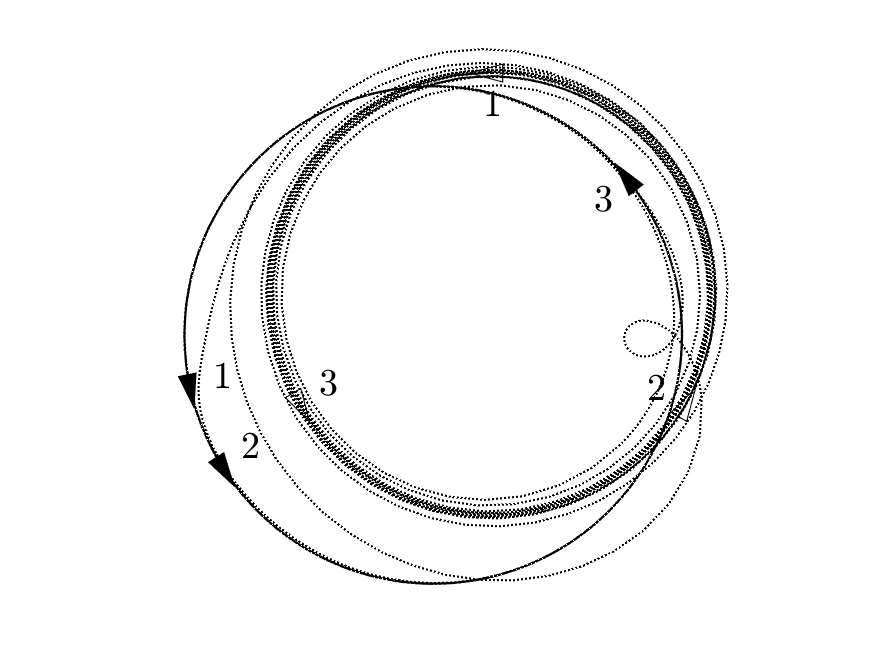}
    \caption{Three unicycles diverge from the nearby irregular formation and eventually converge to a regular one.}
    \label{fig-3vehicles-irregular}
\end{figure}

An illustrative example is provided in~Fig.~\ref{fig-3vehicles}.
While irregular formations are not analyzed, simulations suggest that they are unlikely to be stable. 
As shown in~Fig.~\ref{fig-3vehicles-irregular}, a perturbation cause the group to diverge from the nearby irregular formation and eventually converge to a regular one. 
\begin{remark}
Further research is warranted for the stability analysis of such spontaneous circular formation with larger agent groups. For different $n\geq3$ and for bearing configurations associated with irregular formations, extensive numerical computations invariably yield eigenvalues with positive real parts, regardless of how the gain parameter $k$ is tuned to make $\tr{\hat{A}_R}$ sufficiently negative. This finding leads us to conjecture that only the formation of regular formations is locally asymptotically stable, whereas the formation of irregular ones is unstable. 
\end{remark}

\subsection{A Special Case of \texorpdfstring{$n>3$}{n>3}}
In this part, a special case is studied, where all unicycles are equally spaced in cyclic order, and the following theorem is established. 
\begin{theorem}\label{thm-negeigen}
    For $n>3$, $k=\pm2v$, and $\bar{\xi}\in\mathcal{C}_\pi^-$ satisfying $|\bar{\alpha}_i|=\pi/n$ for all $i$, there exists sufficiently small $\epsilon$ such that the neighborhood
    \begin{equation*}
        \mathcal{U}_\epsilon(\bar{\xi})\coloneqq\left\{\bar{\xi}^\prime\in\mathcal{C}_\pi^- : \|\bar{\xi}^\prime-\bar{\xi}\|<\epsilon\right\}
    \end{equation*} 
    is a locally asymptotically stable invariant set, i.e., a perturbed group either converges to the original regular and equally spaced counterclockwise (resp. clockwise) formation or, in most cases, converges to a regular and `nearly’ equally spaced counterclockwise (resp. clockwise) formation when $k=-2v$ (resp. $k=2v$).
\end{theorem}
A proof sketch is placed in Appendix~\ref{app-negeigen}.

\section{Almost global convergence case study}\label{sec-global}
In this section we propose an alternative control law to achieve almost global convergence for $n=2$. Another appealing advantage is that the bearing spacing between the two unicycles at equilibrium depends solely on their initial heading difference. \\
Consider the following variant of~\eqref{eq-omega-1}:
\begin{equation*}\label{eq-omega-2}
    \omega_i = \frac{v}{\rho_i}\left(\sin{\alpha_i}+\sin{(\alpha_i+\beta_i)}\right)+k\left(\cos{\alpha_i}+\cos{(\alpha_i+\beta_i)}\right).
\end{equation*}
The analysis is simplified since $\rho_1=\rho_2$, $\alpha_2=\alpha_1+\beta_1$, and $\alpha_1=\alpha_2+\beta_2$. Inserting the control law into \eqref{eq-relative-kinematic-1} gives
\begin{subequations}\label{eq-relative-kinematic-n=2}
    \begin{align}
        \dot{\rho}_1 &= -v \left(\cos{\alpha_1}+\cos{(\alpha_1+\beta_1)}\right), \label{eq-relative-kinematic-n=2-rho} \\
        \dot{\alpha}_1 &=  -k\left(\cos{\alpha_1}+\cos{(\alpha_1+\beta_1)}\right), \label{eq-relative-kinematic-n=2-alpha} \\
        \dot{\beta}_1 &= 0, \\
        \dot{\rho}_2 &= -v \left(\cos{\alpha_2}+\cos{(\alpha_2+\beta_2)}\right), \notag \\
        \dot{\alpha}_2 &=  -k\left(\cos{\alpha_2}+\cos{(\alpha_2+\beta_2)}\right), \notag  \\
        \dot{\beta}_2 &= 0. \notag 
    \end{align}
\end{subequations}
Since the motion equations are decoupled, we can drop the indices to simplify notation and proceed with~\eqref{eq-relative-kinematic-n=2}. $\beta(0)=-\pi$ implies that the unicycles initially have identical heading, and~\eqref{eq-relative-kinematic-n=2-rho}--\eqref{eq-relative-kinematic-n=2-alpha} reduce to $\dot{\rho}=0$ and $\dot{\alpha}=0$. Then, depending on the value of $\alpha(0)$, the unicycles advance either on a line or in parallel. We establish the following results whenever $\beta(0)\neq-\pi$. 
\begin{theorem}\label{thm-global-stable}
    Consider the two-unicycle system~\eqref{eq-relative-kinematic-n=2}. Suppose $\theta_1(0)
    \neq\theta_2(0)\Leftrightarrow \beta(0)\neq-\pi$, then the unicycles asymptotically converge to a circular formation. The resulting rotation is counterclockwise (resp. clockwise) if $k<0$ (resp. $k>0$). The bearing spacing at equilibrium is
    \begin{equation}\label{eq-equilibrium-n=2}
        \bar{\alpha}=\left\{
            \begin{aligned}
        &\frac{\pi-\beta(0)}{2}\in(0,\pi), ~\text{if}~k<0; \\
        &\frac{-\pi-\beta(0)}{2}\in(-\pi,0), ~\text{if}~k>0.
        \end{aligned}\right.
    \end{equation}
\end{theorem}
\proof Since $\dot{\beta}=0$, we have $\beta\equiv\beta(0)$. Let $(\bar{\rho}, \bar{\alpha}, \beta(0))$ be an equilibrium to~\eqref{eq-relative-kinematic-n=2}. Then, it follows that $\bar{\alpha}+(\bar{\alpha}+\beta(0))=\pi~\text{mod}~2\pi\Rightarrow\bar{\alpha}=(\pi-\beta(0))/2~\text{mod}~\pi$. Since $\beta(0)\neq-\pi$, after wrapping $\alpha$ into the range $[-\pi,\pi)$, we obtain
\begin{equation*}
    \bar{\alpha} \in 
    \left\{
        \frac{\pi-\beta(0)}{2},\;
        \frac{-\pi-\beta(0)}{2}
    \right\}.
\end{equation*}
We show that~\eqref{eq-equilibrium-n=2} is an almost globally asymptotically stable equilibrium to~\eqref{eq-relative-kinematic-n=2-alpha}. Consider the following Lyapunov function with $\beta(0)\neq-\pi$:
\begin{align*}
    V(\alpha) &= \frac{1}{2}\left(\cos{\alpha}+\cos{(\alpha+\beta(0))}\right)^2 \\
    &~\quad + \frac{1}{2}\left(\sin{\alpha}+\sin{(\alpha+\beta(0))}-2\sin{\bar{\alpha}}\right)^2.
\end{align*}
We have $V(\alpha)\geq0$ and $V(\alpha)=0$ if and only if $\alpha=\bar{\alpha}$ (given that $\alpha\in[-\pi,\pi)$). Its time derivative satisfies 
$\dot{V}(\alpha) = 2k\sin{\bar{\alpha}}\left(\cos{\alpha}+\cos{(\alpha+\beta(0))}\right)^2\leq0$.
Since the invariant set such that $\dot{V}(\alpha)=0$ is the singleton~\eqref{eq-equilibrium-n=2}, LaSalle's theorem confirms that~\eqref{eq-equilibrium-n=2} is almost globally asymptotically stable. Then, as $t\rightarrow\infty$, $\dot{\rho}\rightarrow0$, so $\rho\rightarrow\bar{\rho}$, which can be any positive constant. Thus, $(\bar{\rho}, \bar{\alpha}, \beta(0))$ is locally asymptotically stable. We conclude by~Proposition~\ref{thm-exist.-condt.} that $(\bar{\rho}, \bar{\alpha}, \beta(0))$ corresponds to a particular circular formation.
\endproof

\section{Extension to distance-dependent neighborhoods}\label{sec-dist}
\begin{figure*}
    \centering
    \includegraphics[width=.95\linewidth]{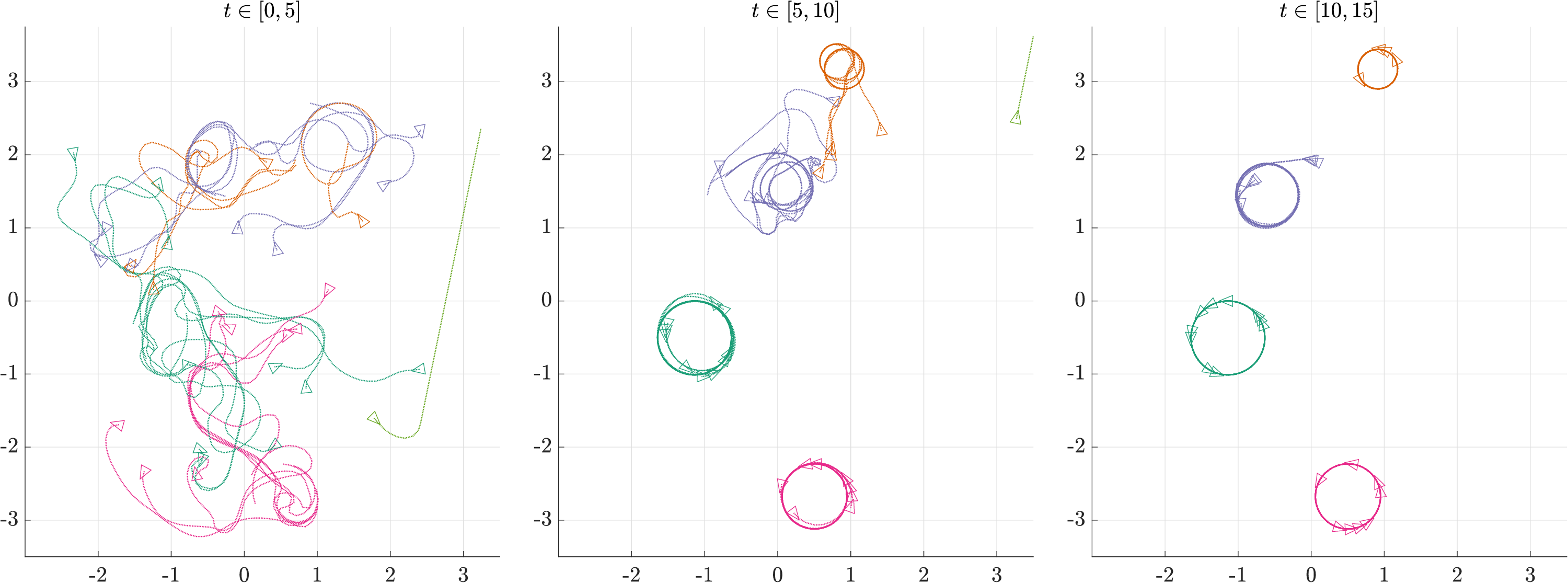}
    \vspace{1em}
    \includegraphics[width=.95\linewidth]{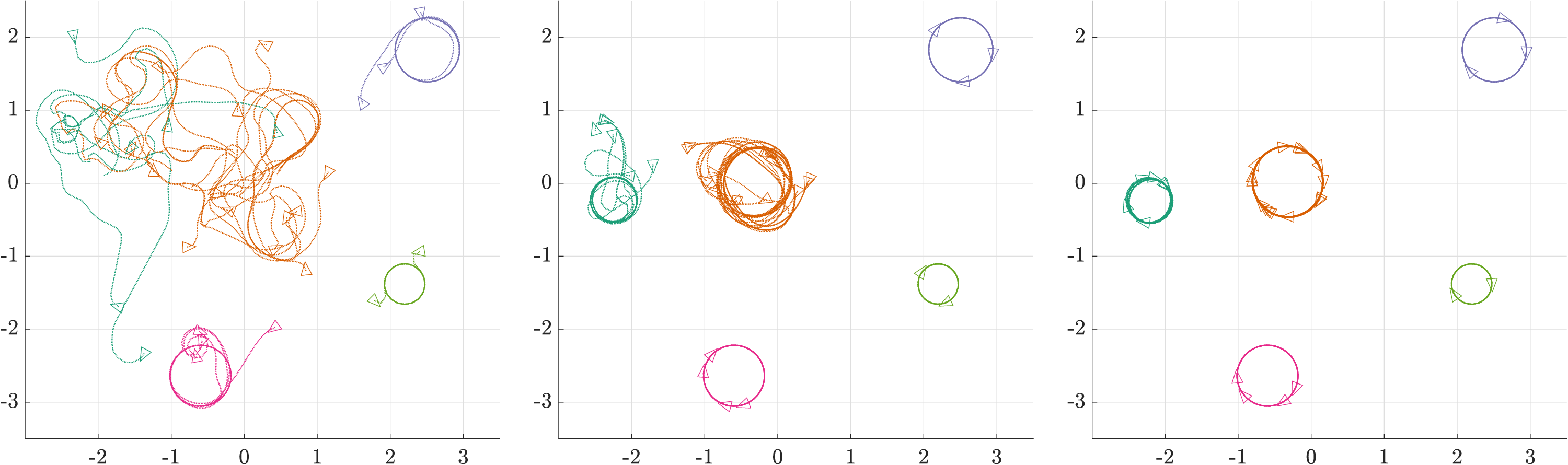}    
    \caption{Due to the limited sensing range $d=1$, $n=30$ unicycles controlled by~\eqref{eq-omega-3} with $v=1$ and $k=\pm5$ converge to several clusters, each forming a circular formation. Under identical initial conditions, $k<0$ results in counterclockwise rotation (top), whereas $k>0$ results in clockwise rotation (below). Almost all clusters have a complete-graph topology; the green cluster in the top picture does not, but it is nearly complete, with only two agents having 8 neighbors and the rest 9.}
    \label{fig-30vehicles}
\end{figure*}

In this section we extend the control law~\eqref{eq-omega-1} to distance-dependent neighborhoods. Numerical simulations show that the group converges to different clusters, most exhibiting a complete-graph topology. Motivated by this observation, a preliminary stability analysis is conducted for the case of a complete graph with $n=3$.
\\
In cyclic pursuit, for each $i\in\mathcal{I}$ we let $\rho_i$ and $\alpha_i$ be the distance and bearing to the next one, respectively. Define $\rho_{ij}$ and $\alpha_{ij}$ in a similar manner for all $j\in\mathcal{N}_i(t)\coloneqq\left\{j\in\mathcal{I}\setminus\{i\}|\rho_{ij}\leq d \right\}$, where $d>0$ is the common sensing range. An extended control law is then given by
\begin{equation}\label{eq-omega-3}
    \omega_i = \frac{1}{\left|\mathcal{N}_i(t)\right|}\sum_{j\in\mathcal{N}_i(t)}{\frac{2v}{\rho_{ij}}\sin{\alpha_{ij}}+\frac{k}{\rho_{ij}}\left(\cos{\alpha_{ij}}+\cos{\alpha_{ji}}\right)},
\end{equation}
which is applied whenever $\mathcal{N}_i(t)\neq\emptyset$; otherwise, we let $\omega_i=0$. 
\\
In order to establish an equilibrium analysis with~\eqref{eq-omega-3}, without loss of generality, we drop the time dependence of $\mathcal{N}_i(t)$ and assume a static topology. Then the new motion equations are analogously given by
\begin{subequations}\label{eq-relative-kinematic-dist}
\begin{align}
    \dot{\rho}_{ij} &= -v\left(\cos{\alpha_{ij}}+\cos{\alpha_{ji}}\right),~i\in\mathcal{I},~j\in\mathcal{N}_i, \\
    \dot{\alpha_{ij}} &= \frac{v}{\rho_{ij}}\left(\sin{\alpha_{ij}}+\sin{\alpha_{ji}}\right)-\omega_i.
\end{align}
\end{subequations}
Define $\mathcal{E}\coloneqq\bigcup_{i\in\mathcal{I}}\left\{(i,j)|j\in\mathcal{N}_i\right\}$ and assume that $\mathcal{G}=\left(\mathcal{I},\mathcal{E}\right)$ is undirected and connected, so $\mathcal{N}_i\neq\emptyset,~\forall i\in\mathcal{I}$. Write $n_\mathsf{e}\coloneqq|\mathcal{E}|/2$, where $|\mathcal{E}|$ is the cardinality of $\mathcal{E}$ and $n_\mathsf{e}$ the number of (bidirectional) edges. Each contributes three coordinates: $\rho_{ij},\alpha_{ij},\alpha_{ji}$, giving a total of $3n_\mathsf{e}$ coordinates. The duplicate coordinate $\rho_{ji}$ does not count since $\rho_{ij}\equiv\rho_{ji}$.
\\
In the rest of the section, we retain the notations from cyclic pursuit, with their definitions understood from context. View \eqref{eq-relative-kinematic-dist} as $\dot{\xi}=\hat{f}(\xi)$, which is lifted to $\mathbb{R}^{3n_\mathsf{e}}$ ($n_\mathsf{e}> n$ in general). In a manner similar to Lemma~\ref{thm-f-iinvariance}, it can be proved that, for arbitrary topologies, there exists a $(3n-3)$-dimensional submanifold $\mathcal{M}\subset\mathbb{R}^{3n_\mathsf{e}}$ that is invariant under $\hat{f}$. The rationale is simple: $3n-3$ coordinates are enough to describe an $n$-unicycle system in relative coordinates. 
\begin{theorem}\label{thm-equilibria-dist}
    An invariant set to the system $\dot{\xi}=\hat{f}(\xi)$ is given by
    \begin{equation}\begin{aligned}
    \Xi \coloneqq \big\{&\bar{\xi}\in\mathcal{M} : \bar{\alpha}_{ij}+\bar{\alpha}_{ji}=\pi~\text{mod}~2\pi~\text{and} \\ 
    &\exists\bar{s}~\text{s.t.}~{\sin{\bar{\alpha}}_{ij}}/{\bar{\rho}_{ij}}=\bar{s},~\forall i\in\mathcal{I},\forall j\in\mathcal{N}_i
    \big\}. \label{eq-set-equilibria-dist}
    \end{aligned}\end{equation}
    A circular formation is realized whenever $\bar{s} \neq 0$, counterclockwise for $\bar{s}>0$ and clockwise for $\bar{s}<0$. 
\end{theorem} 
The relevant results in cyclic pursuit naturally generalize to connected graphs (and indeed to strongly connected graphs), and thus the proof is omitted. 

Other feasible invariant sets might exist; however, they are not of interest.
\\
For $\bar{\xi}\in\Xi$, write $\tilde{\xi}=\xi-\bar{\xi}$ and denote the linearized system as $\dot{\tilde{\xi}}=\hat{A}\tilde{\xi}$. The explicit expression of $\hat{A}\in\mathbb{R}^{3n_\mathsf{e}\times 3n_\mathsf{e}}$ is tedious and thus omitted. By arguments similar to those in Theorem~\ref{thm-unstable}, we can show that the formation of collinear formations is not asymptotically stable. The following theorem is thus given without proof.
\begin{theorem}
    The invariant set $\mathcal{L}$ is not asymptotically stable, i.e., the formation of collinear formations is not asymptotically stable.
\end{theorem} 
The complement $\mathcal{C}=\Xi\setminus\mathcal{L}$ comprises two disjoint components, $\mathcal{C}^-$ and $\mathcal{C}^+$; see~\eqref{eq-C-two-parts} for their definitions. The latter is proven to be unstable.
\begin{theorem}\label{thm-dist-unstable}
    The invariant set $\mathcal{C}^+$ is unstable, i.e., counterclockwise (resp. clockwise) formations are unstable when $k>0$ (resp. $k<0$).
\end{theorem}
\proof
We show that the matrix $\hat{A}$ has a positive trace for every $\bar{\xi}\in\mathcal{C}^+$. Write $\bar{c}_{ij}\coloneqq \cos{\bar{\alpha}_{ij}}/\bar{\rho}_{ij}$, the trace of $\hat{A}$ satisfies
\begin{equation*}\begin{aligned}
    \tr{\hat{A}} &= \sum_{i\in\mathcal{I}}\sum_{j\in\mathcal{N}_i}{ \left. \frac{\partial\dot{\rho}_{ij}}{\partial\rho_{ij}} \right|_{\bar{\xi}} + \left. \frac{\partial\dot{\alpha}_{ij}}{\partial\alpha_{ij}} \right|_{\bar{\xi}}} \\
    &=\sum_{i\in\mathcal{I}}\sum_{j\in\mathcal{N}_i}\left( v\bar{c}_{ij} -\frac{1}{\left|\mathcal{N}_i\right|}\sum_{j^\prime\in\mathcal{N}_i}\left(2v\bar{c}_{ij^\prime}-k\bar{s}\right) \right) \\
    &= \sum_{i\in\mathcal{I}}\sum_{j\in\mathcal{N}_i}\left(v\bar{c}_{ij} + k\bar{s}\right) -\sum_{i\in\mathcal{I}}\sum_{j\in\mathcal{N}_i}\left(\frac{2v}{\left|\mathcal{N}_i\right|}\sum_{j^\prime\in\mathcal{N}_i}\bar{c}_{ij^\prime}\right).
\end{aligned}\end{equation*}
Note that $\frac{\partial\dot{\rho}_{ij}}{\partial\rho_{ij}}\equiv0$. It holds that
\begin{equation*}\begin{aligned} 
    &\sum_{j\in\mathcal{N}_i}\left(\frac{1}{\left|\mathcal{N}_i\right|}\sum_{j^\prime\in\mathcal{N}_i}\bar{c}_{ij^\prime}\right) = \sum_{j\in\mathcal{N}_i}\frac{1}{\left|\mathcal{N}_i\right|}\cdot\sum_{j^\prime\in\mathcal{N}_i}\bar{c}_{ij^\prime} = \sum_{j\in\mathcal{N}_i}\bar{c}_{ij}.
\end{aligned}\end{equation*}
Since the graph is undirected and $\bar{c}_{ij}+\bar{c}_{ji}=0$, we yield $\sum_{i\in\mathcal{I}}\sum_{j\in\mathcal{N}_i}{\bar{c}_{ij}}=0$ and  $\tr{\hat{A}}=\sum_{i\in\mathcal{I}}\sum_{j\in\mathcal{N}_i}k\bar{s}$. As a result, $\tr{\hat{A}}>0,~\forall\bar{\xi}\in\mathcal{C}^+$.
\endproof

Regarding the stability of $\mathcal{C}^-$, it is perhaps unnecessary to conduct a general analysis under arbitrary topologies. As suggested by numerical simulations, in most cases, the resulting clusters possess a complete-graph topology; see Fig.~\ref{fig-30vehicles}. The stability in this case is thus of particular interest. On the other hand, $n_\mathsf{e}$ in this case attains its maximal value, and the system $\dot{\xi}=\hat{f}(\xi)$ is lifted to its highest possible dimension. \\
A complete graph when $n=3$ is a bidirectional cyclic pursuit. For this case, we can instead adopt the coordinate representation in cyclic pursuit: 
\begin{equation*}
    \rho_{i(i+1)} = \rho_i,~\alpha_{i(i+1)} = \alpha_i,~ \alpha_{(i+1)i} = \gamma_i \coloneqq \alpha_i+\beta_i. 
\end{equation*}
In this way, the established analysis approach can be carried over to complete graph with $n=3$. For brevity, the routine intermediate derivations arising in this process are omitted. Now, for any $\bar{\xi}\in\mathcal{C}^-$ (expressed in the new coordinate representation), the expression of $\hat{A}$ can be derived as
\begin{align*}
    \hat{A} &= \begin{bmatrix}
        A_1 & B_2 & C_3 \\
        C_1 & A_2 & B_3 \\
        B_1 & C_2 & A_3
    \end{bmatrix},~ 
    A_i = \begin{bmatrix}
        0 & v\bar{s}\rho_i & v\bar{s}\rho_i \\
        -v\bar{s}/\bar{\rho}_i & \bar{k}\bar{s} & \bar{k}\bar{s}-v\bar{c}_i \\
        -v\bar{s}/\bar{\rho}_i & \bar{k}\bar{s}+v\bar{c}_i & \bar{k}\bar{s}
    \end{bmatrix}, \\
    B_i &= \begin{bmatrix}
        0 & 0 & 0 \\ 0 & 0 & 0 \\
        v\bar{s}/\rho_i & \bar{k}\bar{s}-v\bar{c}_i & \bar{k}\bar{s}
    \end{bmatrix},~ 
    C_i = \begin{bmatrix}
        0 & 0 & 0 \\ 
        v\bar{s}/\bar{\rho}_i & \bar{k}\bar{s} & \bar{k}\bar{s}+v\bar{c}_i \\
        0 & 0 & 0
    \end{bmatrix}, \\
    \bar{k} &\coloneqq k/2, ~\bar{c}_i \coloneqq \cos{\bar{\alpha}_i}/\bar{\rho}_i, ~i\in\{1,2,3\}.
\end{align*}
The following corollary, following the same rationale and reasoning as the corresponding results in cyclic pursuit, is presented without proof.
\begin{corollary}\label{thm-criterion-comple}
    (i) For all $\bar{\xi}\in\mathcal{C}^-$, the matrix $\hat{A}$ has a pair of purely imaginary eigenvalues $\pm j2v\bar{s}$ and four zero eigenvalues.
    (ii) If, for all $\bar{\xi}\in\mathcal{C}_{\pi}^-$ (resp. $\bar{\xi}\in\mathcal{C}_{2\pi}^-$), the remaining three eigenvalues have negative real parts, then $\mathcal{C}_\pi^-$ (resp. $\mathcal{C}_{2\pi}^-$) is locally asymptotically stable.
\end{corollary}
\begin{remark}
    The invariant set $\mathcal{C}^-$ decomposes into two maximal connected subsets, $\mathcal{C}_\pi^-$ and $\mathcal{C}_{2\pi}^-$; see~\eqref{eq-C-of-pi} for their definitions. Since the agents are indistinguishable in the complete graph, the stability is independent of their arrangements along the circle; the two subsets necessarily share the same type of stability: If one is stable (resp. unstable), then so is the other. Nevertheless, it is necessary to treat them separately since they are disjoint, and intuitively an equilibrium lying close to $\mathcal{C}_\pi^-$ is unlikely to converge to $\mathcal{C}_{2\pi}^-$ (and vice versa) if both are stable.
\end{remark}
Then, we transform the matrix $\hat{A}$ into the form~\eqref{eq-A-lower-form} to reduce the model.
\begin{corollary}
    For any $\bar{\xi}\in\mathcal{C}^-$, the matrix $\hat{A}$ can be transformed into the form~\eqref{eq-A-lower-form} via a similarity transformation. Let $v=1$ without loss of generality, then  $\hat{A}_R$ is given by
    \begin{align*}
        \hat{A}_R &= \begin{bmatrix}
            \bar{A}_1 & \bar{B} & \bar{C}_3 \\
            \bar{C}_1 & \bar{A}_2 & \bar{B} \\
            \bar{B} & \bar{C}_2 & \bar{A}_3
        \end{bmatrix}, ~\bar{B} = |\bar{s}|\begin{bmatrix}
            0 & 0 \\ 1 & \bar{k}
        \end{bmatrix}, \\     
        \bar{A}_i &= |\bar{s}|\begin{bmatrix}
            \cot{|\bar{\alpha}_i|} & \cot^2{|\bar{\alpha}_i|}-\bar{k}\cot{|\bar{\alpha}_i|}+1 \\
            -2 & 2\bar{k}-\cot{|\bar{\alpha}_i|}
        \end{bmatrix}, \\
        \bar{C}_i &= |\bar{s}|\begin{bmatrix}
            -\cot{|\bar{\alpha}_{i+1}|} & -\cot{|\bar{\alpha}_{i+1}|}\left(\bar{k}+\cot{|\bar{\alpha}_i|}\right) \\
            1 & \bar{k}+\cot{|\bar{\alpha}_i|}
        \end{bmatrix}, \\
        \bar{k} &< 0, ~|\bar{\alpha}_i|\in(0,\pi), ~i\in\{1,2,3\}.
    \end{align*}
\end{corollary}
\proof The proof parallels that of Lemma~\ref{thm-A-reduced}, except that the similarity transformation now uses a modified matrix
$T=\blkdiag{\left(T_1,T_2,T_3\right)}$, where
\begin{equation*}
    T_i= \begin{bmatrix}
        {1}/{\bar{\rho}_i} & -\cot{\alpha}_i & 0 \\ 0 & 1 & 0 \\ 0 & 1 & 1
    \end{bmatrix}, ~i\in\{1,2,3\}.
\end{equation*}
This completes the proof.
\endproof

We then derive the characteristic polynomial of $\hat{A}_R/|\bar{s}|$ through laborious computations; the procedure is omitted due to its cumbersome nature. This leads to the final theorem of the paper.
\begin{theorem}\label{thm-dist-n=3}
    Consider the system~\eqref{eq-relative-kinematic-dist} with a complete graph. For $n=3$, both invariant sets $\mathcal{C}_\pi^-$ and $\mathcal{C}_{2\pi}^-$ are locally asymptotically stable, i.e., the formation of counterclockwise (resp. clockwise) formations is locally asymptotically stable when $k<0$ (resp. $k>0$).
\end{theorem}
\proof
The characteristic polynomial of $\hat{A}_R/|\bar{s}|$ has the form
\begin{equation*}\begin{aligned}
    P_{\bar{\xi}}(\lambda) &= \lambda^6+b_5\lambda^5+b_4\lambda^4+b_3\lambda^3 + b_2\lambda^2+b_1\lambda \\
    &= \left(\lambda^3+4\lambda\right)\left(\lambda^3+a_2\lambda^2+a_1\lambda+a_0\right),
\end{aligned}\end{equation*}
and the expression of $b_5,b_4,b_3$, through laborious computations, are derived as $b_5=-6\bar{k}$,
\begin{align*}
    b_4 &= 9{\bar{k}}^2 + 6 +  \sum_{i=1}^3{\cot^2{|\bar{\alpha}_i|}} + \sum_{i=1}^3{\cot{|\bar{\alpha}_i|}\cot{|\bar{\alpha}_{i+1}|}}, \\
    b_3 &= -4{\bar{k}}^3 -4\bar{k}\sum_{i=1}^3{\cot^2{|\bar{\alpha}_i|}} 
    \\ &~\quad 
    -6\bar{k}\sum_{i=1}^3{\cot{|\bar{\alpha}_i|}\cot{|\bar{\alpha}_{i+1}|}} - 30\bar{k},
\end{align*}
where $\sum_{i=1}^3{\cot{|\bar{\alpha}_i|}\cot{|\bar{\alpha}_{i+1}|}}=1$ since $\sum_{i=1}^3{|\bar{\alpha}_i|}=\pi$ or $2\pi$. Notice that
\begin{equation*}
    a_2 = b_5 = -6\bar{k},~ a_1 = b_4-4, ~ a_0 = b_3-4b_5.
\end{equation*}
Write $c^2\coloneqq 3+\sum_{i=1}^3{\cot^2{|\bar{\alpha}_i|}}$, then $a_1=c^2+9{\bar{k}}^2$ and $a_0 = -4\bar{k}c^2-4{\bar{k}}^3$.
With this particular representation, we can see that
\begin{align*}
    \lambda^3+a_2\lambda^2+a_1\lambda+a_0 = \left(\lambda-4\bar{k}\right)\left(\left(\lambda-\bar{k}\right)^2+c^2\right).
\end{align*}
Hence, the remaining three eigenvalues of $\hat{A}$ have negative real parts for all $\bar{\xi}\in\mathcal{C}^-$. The rest of the theorem then follows from Corollary~\ref{thm-criterion-comple}.
\endproof

\section{Concluding remarks}\label{sec-conclusion}
To conclude, we have investigated the emergence of circular formation for unicycles in cyclic pursuit. A local control law was proposed to achieve spontaneous circular formation. The closed-loop system possesses infinitely many non-isolated equilibria, which are not asymptotically stable and form a disconnected invariant set. Consequently, a stability analysis was developed at the level of maximal connected invariant sets rather than individual equilibria. Based on the derived local stability criterion, we established a set of stability results for small agent groups ($n\leq3$), providing insights into the stability property of such spontaneous circular formation. 
\\
In addition, we have provided two extensions of the proposed control law. First, a modification was introduced to achieve almost global convergence for $n=2$. Second, the control law was extended to distance-dependent neighborhoods. The topology is rendered undirected, and we proved that counterclockwise (resp. clockwise) formations are unstable whenever $k>0$ (resp. $k<0$). This phenomenon is also observed in the simulations of cyclic pursuit, but cannot be fully analyzed there. The group converges to several clusters, most exhibiting a complete-graph topology. Thus, we developed a preliminary stability analysis for the complete-graph case with $n=3$.
\\
The paper opens up various interesting problems at both microscopic and macroscopic scales. At the microscopic scale, it warrants further investigation whether a more scalable approach can be developed to analyze stability with larger agent groups. While attaining a fully exact understanding of such complex collective behavior is expected to be highly nontrivial, it is perhaps not the most immediate priority. Indeed, relatively accessible problems arise at the macroscopic scale, such as investigating whether the observed collective behaviors exhibit convergence in some statistical sense. Addressing these problems sheds light on the mechanisms underlying complex natural phenomena such as fish milling and starling murmuration.

\bibliographystyle{agsm}
\bibliography{refs}

@article{marshall2004formations,
  title={Formations of vehicles in cyclic pursuit},
  author={Marshall, Joshua~A. and Broucke, Mireille~E. and Francis, Bruce~A.},
  journal={IEEE Transactions on Automatic Control},
  volume={49},
  number={11},
  pages={1963--1974},
  year={2004},
  publisher={IEEE}
}

@article{marshall2006pursuit,
  title={Pursuit formations of unicycles},
  author={Marshall, Joshua~A. and Broucke, Mireille~E. and Francis, Bruce~A.},
  journal={Automatica},
  volume={42},
  number={1},
  pages={3--12},
  year={2006},
  publisher={Elsevier}
}

@article{zheng2009ring,
  title={Ring-coupled unicycles: Boundedness, convergence, and control},
  author={Zheng, Ronghao and Lin, Zhiyun and Yan, Gangfeng},
  journal={Automatica},
  volume={45},
  number={11},
  pages={2699--2706},
  year={2009},
  publisher={Elsevier}
}

@article{zheng2015distributed,
  title={Distributed control for uniform circumnavigation of ring-coupled unicycles},
  author={Zheng, Ronghao and Lin, Zhiyun and Fu, Minyue and Sun, Dong},
  journal={Automatica},
  volume={53},
  pages={23--29},
  year={2015},
  publisher={Elsevier}
}

@article{shi2021distributed,
  title={Distributed circumnavigation control of autonomous underwater vehicles based on local information},
  author={Shi, Linlin and Zheng, Ronghao and Liu, Meiqin and Zhang, Senlin},
  journal={Systems \& Control Letters},
  volume={148},
  pages={104873},
  year={2021},
  publisher={Elsevier}
}

@article{el2012distributed,
  title={Distributed circular formation stabilization for dynamic unicycles},
  author={El-Hawwary, Mohamed~I. and Maggiore, Manfredi},
  journal={IEEE Transactions on Automatic Control},
  volume={58},
  number={1},
  pages={149--162},
  year={2012},
  publisher={IEEE}
}

@article{yu2022decentralized,
  title={Decentralized circular formation control of nonholonomic mobile robots under a directed sensor graph},
  author={Yu, Xiao and Su, Rong},
  journal={IEEE Transactions on Automatic Control},
  volume={68},
  number={6},
  pages={3656--3663},
  year={2022},
  publisher={IEEE}
}

@article{sepulchre2007stabilization,
  title={Stabilization of planar collective motion: All-to-all communication},
  author={Sepulchre, Rodolphe and Paley, Derek A. and Leonard, Naomi E.},
  journal={IEEE Transactions on Automatic Control},
  volume={52},
  number={5},
  pages={811--824},
  year={2007},
  publisher={IEEE}
}

@article{sepulchre2008stabilization,
  title={Stabilization of planar collective motion with limited communication},
  author={Sepulchre, Rodolphe and Paley, Derek A. and Leonard, Naomi E.},
  journal={IEEE Transactions on Automatic Control},
  volume={53},
  number={3},
  pages={706--719},
  year={2008},
  publisher={IEEE}
}

@article{galloway2018collective,
  title={Collective motion under beacon-referenced cyclic pursuit},
  author={Galloway, Kevin S. and Dey, Biswadip},
  journal={Automatica},
  volume={91},
  pages={17--26},
  year={2018},
  publisher={Elsevier}
}

@article{seyboth2014collective,
  title={Collective circular motion of unicycle type vehicles with nonidentical constant velocities},
  author={Seyboth, Georg S. and Wu, Jingbo and Qin, Jiahu and Yu, Changbin and Allg{\"o}wer, Frank},
  journal={IEEE Transactions on Control of Network Systems},
  volume={1},
  number={2},
  pages={167--176},
  year={2014},
  publisher={IEEE}
}

@ARTICLE{xie2025circular,
  author={Xie, Jingjing and Lan, Weiyao and Tong, Feng and Yu, Xiao},
  journal={IEEE Transactions on Control of Network Systems}, 
  title={Circular Formation Control for Networked Dynamic Autonomous Underwater Vehicles: Estimator-Based Approaches}, 
  year={2025},
  volume={12},
  number={3},
  pages={2016-2028}
}

@article{tripathy2024convergence,
  title={On convergence results for nonlinear cyclic pursuit strategies},
  author={Tripathy, Twinkle and Shima, Tal},
  journal={Automatica},
  volume={159},
  pages={111315},
  year={2024},
  publisher={Elsevier}
}

@book{marshall2005coordinated,
  title={Coordinated autonomy: Pursuit formations of multivehicle systems},
  author={Marshall, Joshua Alexander},
  year={2005},
  publisher={University of Toronto}
}

@article{li2022differential,
  title={A differential game approach to intrinsic formation control},
  author={Li, Yibei and Hu, Xiaoming},
  journal={Automatica},
  volume={136},
  pages={110077},
  year={2022},
  publisher={Elsevier}
}

@article{zhang2020intrinsic,
  title={An intrinsic approach to formation control of regular polyhedra for reduced attitudes},
  author={Zhang, Silun and He, Fenghua and Hong, Yiguang and Hu, Xiaoming},
  journal={Automatica},
  volume={111},
  pages={108619},
  year={2020},
  publisher={Elsevier}
}

@article{sen2019circumnavigation,
  title={Circumnavigation on multiple circles around a nonstationary target with desired angular spacing},
  author={Sen, Arijit and Sahoo, Soumya Ranjan and Kothari, Mangal},
  journal={IEEE Transactions on Cybernetics},
  volume={51},
  number={1},
  pages={222--232},
  year={2019},
  publisher={IEEE}
}

@article{yu2018cooperative,
  title={Cooperative moving-target enclosing of networked vehicles with constant linear velocities},
  author={Yu, Xiao and Ding, Ning and Zhang, Aidong and Qian, Huihuan},
  journal={IEEE Transactions on Cybernetics},
  volume={50},
  number={2},
  pages={798--809},
  year={2018},
  publisher={IEEE}
}

@article{chun2020multi,
  title={Multi-targets localization and elliptical circumnavigation by multi-agents using bearing-only measurements in two-dimensional space},
  author={Chun, Shaoheng and Tian, Yu-Ping},
  journal={International Journal of Robust and Nonlinear Control},
  volume={30},
  number={8},
  pages={3250--3268},
  year={2020},
  publisher={Wiley Online Library}
}

@article{zhang2025angle,
  author = {Zhang, Chunyan and Chen, Ye and Liang, Yuan and Kong, Deren and Li, Yinya and Sheng, Andong},
  title = {Angle Rigidity-Based Distributed Circular Formation Control in 2D Space},
  journal = {International Journal of Robust and Nonlinear Control},
  volume = {35},
  number = {12},
  pages = {5077-5093},
  year = {2025},
  publisher={Wiley Online Library}
}

@article{sun2018circular,
  title={Circular formation control of multiple unicycle-type agents with nonidentical constant speeds},
  author={Sun, Zhiyong and De Marina, Hector Garcia and Seyboth, Georg S and Anderson, Brian DO and Yu, Changbin},
  journal={IEEE Transactions on Control Systems Technology},
  volume={27},
  number={1},
  pages={192--205},
  year={2018},
  publisher={IEEE}
}

@article{dou2020target,
  title={Target localization and enclosing control for networked mobile agents with bearing measurements},
  author={Dou, Liya and Song, Cheng and Wang, Xiaofan and Liu, Lu and Feng, Gang},
  journal={Automatica},
  volume={118},
  pages={109022},
  year={2020},
  publisher={Elsevier}
}

@article{miao2017cooperative,
  title={Cooperative circumnavigation of a moving target with multiple nonholonomic robots using backstepping design},
  author={Miao, Zhiqiang and Wang, Yaonan and Fierro, Rafael},
  journal={Systems \& Control Letters},
  volume={103},
  pages={58--65},
  year={2017},
  publisher={Elsevier}
}

@article{zhang2020distributed,
  title={Distributed finite-time control for coordinated circumnavigation with multiple agents under directed topology},
  author={Zhang, Chunyan and Li, Yinya and Qi, Guoqing and Sheng, Andong},
  journal={Journal of the Franklin Institute},
  volume={357},
  number={16},
  pages={11710--11729},
  year={2020},
  publisher={Elsevier}
}

\setcounter{section}{0}
\setcounter{equation}{0}
\renewcommand{\thesection}{Appendix} 
\renewcommand{\thesubsection}{\Alph{section}.\arabic{subsection}}
\renewcommand{\theequation}{\Alph{section}.\arabic{equation}}

\section{Deferred proofs}

\subsection{Proof of Proposition~\ref{thm-exist.-condt.-global}}\label{app-exist.-condt.-global}
\emph{Sufficiency:} Let $\mathtt{v}_j=\left[\cos{\theta_j},\sin{\theta_j}\right]^\top$ be the unit heading vector. A counterclockwise (resp. clockwise) rotation of $\mathtt{v}_j$ by $\pi/2$ yields $\mathtt{r}_j$ (resp. $-\mathtt{r}_j$). Since $\|\tilde{z}_j\|$ and $\|\tilde{\mathtt{r}}_j\|$ are nonzero,  $\tilde{z}_j\times\tilde{\mathtt{r}}_j=0\Leftrightarrow \tilde{z}_j=r_j\tilde{\mathtt{r}}_j\Leftrightarrow z_j-r_j\mathtt{r}_j=z_1-r_j\mathtt{r}_1$, where $r_j={\tilde{z}_j^\top\tilde{\mathtt{r}}_j}/{\|\tilde{\mathtt{r}}_j\|^2}\neq0$. Thus, unicycles $1$ and $j\geq2$ travel along a circle, counterclockwise (resp. clockwise) for $r_j<0$ (resp. $r_j>0$), as illustrated in Fig.~\ref{fig-product-zero}. Since $r\coloneqq r_2=r_j$ for all $j\geq3$, we conclude that all vehicles travel along a common circle at angular speed $v/r$, counterclockwise (resp. clockwise) for $r<0$ (resp. $r>0$). This completes the sufficiency. \\
\emph{Necessity:} This part is straightforward and thus omitted. \\
We next prove the final statement of the proposition. For ease of reference, we introduce additional notations. Write $q=\left[q_1^\top,q_2^\top,\ldots,q_n^\top\right]^\top$ with $q_i=\left[z_i^\top,\theta_i\right]^\top$, and 
\begin{align*}
    h_{j,1}(q) &\coloneqq \frac{\tilde{z}_2^\top\tilde{\mathtt{r}}_2}{\|\tilde{\mathtt{r}}_2\|^2} - \frac{\tilde{z}_j^\top\tilde{\mathtt{r}}_j}{\|\tilde{\mathtt{r}}_j\|^2}, ~j\in\{3,4,\ldots,n\}, \\
    h_{j,2}(q) &\coloneqq \tilde{z}_j\times\tilde{\mathtt{r}}_j, ~j\in\{2,3,\ldots,n\}.
\end{align*}
Let $h(q) = \left[h_{2,2}(q),\ldots,h_{n,1}(q),h_{n,2}(q)\right]^\top\in \mathbb{R}^{2n-3}$. The submanifold $\mathcal{H}\coloneqq\left\{q\in\mathbb{R}^{3n}|h(q)=0\right\}\subset\mathbb{R}^{3n}$
has dimension $n+3$ if $\frac{\partial h(q)}{\partial q}$ has full rank, and we prove this by induction. \\
The results holds trivially for the base case $n=2$ since $\|\tilde{z}_j\|$ and $\|\tilde{\mathtt{r}}_j\|$ are nonzero. Assume that $\frac{\partial h(q)}{\partial q}$ has full rank for $n=k\geq2$, and consider the case $n=k+1$. Note that, for all $j\le k$, $h_{j,1}(q)$ and $h_{j,2}(q)$ do not depend on $q_{k+1}$. Thus, the first $2k-3$ rows of $\frac{\partial h(q)}{\partial q}$ have zero entries in the columns corresponding to $q_{k+1}$. Let $h_{k+1}(q) = \left[h_{k+1,1}(q), h_{k+1,2}(q)\right]^\top$. A straightforward calculation yields
\begin{align*}
    \frac{\partial h_{k+1}(q)}{\partial z_{k+1}} &= \begin{bmatrix}  
    \displaystyle 
    \frac{\tilde{\mathtt{r}}_{k+1,x}}{\|\tilde{\mathtt{r}}_{k+1}\|^2} & \displaystyle \frac{\tilde{\mathtt{r}}_{k+1,y}}{\|\tilde{\mathtt{r}}_{k+1}\|^2} \\
    -\tilde{\mathtt{r}}_{k+1,y} & \tilde{\mathtt{r}}_{k+1,x}
    \end{bmatrix},
\end{align*}
where $\tilde{\mathtt{r}}_{k+1,x}$ and  $\tilde{\mathtt{r}}_{k+1,y}$ denote the $x$ and $y$ components of $\tilde{\mathtt{r}}_{k+1}$, respectively. Since $\|\tilde{\mathtt r}_{k+1}\|\neq 0$, we have
\begin{align*}
    \det{\frac{\partial h_{k+1}(q)}{\partial z_{k+1}}} = \frac{\tilde{\mathtt{r}}_{k+1,x}^2+\tilde{\mathtt{r}}_{k+1,y}^2}{\|\tilde{\mathtt{r}}_{k+1}\|^2} = 1\neq0.
\end{align*}
Thus, the two new rows $\frac{\partial h_{k+1}(q)}{\partial q_{k+1}}$ are linearly independent of each other and of the previous $2k-3$ rows, which confirms that $\frac{\partial h(q)}{\partial q}$ has full rank for $n=k+1$. 
\endproof

\subsection{Proof of Proposition~\ref{thm-exist.-condt.}} \label{app-exist.-condt.}
\emph{Sufficiency:} If we normalize the angle variables into the range $[-\pi,\pi)$, then only two cases need to be considered: (i) Both $\alpha_i$ and $\alpha_i+\beta_i$ lie in the range $(0,\pi)$ for all $i$, and \eqref{eq-sum-pi-mod-2pi} reduces to $\alpha_i+(\alpha_i+\beta_i)=\pi$; (ii) Both lie in the range $(-\pi,0)$ for all $i$, and \eqref{eq-sum-pi-mod-2pi} reduces to $\alpha_i+(\alpha_i+\beta_i)=-\pi$. We have $r>0$ in case (i) and $r<0$ in case (ii). 
\\
We prove case (i) only; the proof for case (ii) is similar and omitted. First, denote $\mathtt{v}_i=\left[\cos{\theta_i},\sin{\theta_i}\right]^\top$ as the unit heading vector and obtain $\mathtt{r}_i=\left[-\sin{\theta_i},\cos{\theta_i}\right]^\top$ by rotating $\mathtt{v}_i$ counterclockwise by $\pi/2$. Then, as illustrated in Fig.~\ref{fig-sum-pi}, the extensions of $z_i+\mathtt{r}_i$ and $z_{i+1}+\mathtt{r}_{i+1}$ in the global reference frame intersect at a common point, denoted as $\mathcal{O}$. From \eqref{eq-sum-pi-mod-2pi} we have $\alpha_i=\pi-(\alpha_i+\beta_i)$, and thus $\pi/2-\alpha_i=(\alpha_i+\beta_i)-\pi/2$, which implies that $\angle\mathcal{O}(i+1)i=\angle\mathcal{O}i(i+1)$. Thus, each pair of $i$ and $i+1$ travels along a circle.
On the other hand, since $\alpha_i\in(0,\pi)$, we have $2\sin{\alpha_i}/\rho_i=1/r>0$. Assume that $r_i>0$ is the radius of the circle along which $i$ and $i+1$ travel. When $\alpha_i\neq\pi/2$, due to the trigonometric relation, we have 
\begin{equation}\label{eq-trigonometric-relation}
    \frac{\sin{\left(\pi-2\left(\pi/2-\alpha_i\right)\right)}}{\rho_i} = \frac{\sin{\left(\pi/2-\alpha_i\right)}}{r_i},
\end{equation}
which is equivalent to $2\sin{\alpha_i}/\rho_i=1/r_i$. Thus, $r_i=r$. When $\alpha_i=\pi/2$, the metric relation between diameter and radius follows: 
\begin{equation}\label{eq-metric-relation}
    \rho_i=2r_i~\Leftrightarrow~2/\rho_i=1/r_i.
\end{equation}
However, $2\sin{\alpha_i}/\rho_i=2/\rho_i=1/r$, so we still obtain $r_i=r$. Therefore, $r$ is the radius of the circle along which $i$ and $i+1$ travel. Since $2\sin{\alpha_i}/\rho_i=1/r$ for all $i$, all circles have the same radius. Therefore, all unicycles travel along a common circle at positive angular speed $v/r$; the rotation is counterclockwise. 
\\
Finally, while circular formation is naturally feasible for unicycles in cyclic pursuit, we theoretically demonstrate that~\eqref{eq-sum-pi-mod-2pi}--\eqref{eq-same-ratio} are compatible with the constraints $g(\xi)=0$ with $g(\xi)\coloneqq[g_1(\xi),g_2(\xi),g_3(\xi)]^\top$.
Given $\beta_i$'s satisfying $g_3(\xi)=0$, feasible $\alpha_i$'s satisfying~\eqref{eq-sum-pi-mod-2pi} clearly exist. 
Then we show that $g_1(\xi)=0$ (resp. $g_2(\xi)=0$) are satisfied whenever \eqref{eq-sum-pi-mod-2pi}–\eqref{eq-same-ratio} and $g_3(\xi)=0$ are given.
\\
To see this, we first note from~\eqref{eq-sum-pi-mod-2pi} that 
$\sum_{i=1}{2\alpha_i+\beta_i}=n\pi~\text{mod}~2\pi$.
Subtracting it by~\eqref{eq-constr-3} yields that
$2\sum_{i=1}^n{\alpha_i}-n\pi=n\pi~\text{mod}~2\pi\Rightarrow2\sum_{i=1}^n{\alpha_i}=0~\text{mod}~2\pi$. 
Without loss of generality, assume
$r=1/2\Rightarrow\rho_i=\sin{\alpha_i}$, 
then $g_1(\xi)=0$ is verified through a series of manipulations as follows:
\begin{equation*}\begin{aligned}
    &g_1(\xi) = \rho_1\sin{\alpha_1} + \sum_{j=2}^n{\rho_j\sin{(\alpha_j+(j-1)\pi-\sum_{k=1}^{j-1}{\beta_k})}} \\
    &= \rho_1\sin{\alpha_1} + \sum_{j=2}^{n}\rho_j\sin{(\alpha_j+2\sum_{k=1}^{j-1}{\alpha_k})} \\
    &= \sin{\alpha_1}\sin{\alpha_1} + \sum_{j=2}^{n}{\sin{\alpha_j}\sin{(\alpha_j+2\sum_{k=1}^{j-1}{\alpha_k})}} \\
    &= \sin{\alpha_1}\sin{\alpha_1} + \frac{1}{2}\sum_{j=2}^{n}{(\cos{(2\sum_{k=1}^{j-1}{\alpha_k})}-\cos{(2\sum_{k=1}^{j}{\alpha_k})})} \\
    &= \sin{\alpha_1}\sin{\alpha_1} + \frac{1}{2}\cos{2\alpha_1} - \frac{1}{2}\cos{(2\sum_{i=1}^{n}{\alpha_i})} = 0.
\end{aligned}\end{equation*}
The second ``$=$'' is obtained by taking into account $2\alpha_i+\beta_i=\pi~\text{mod}~2\pi$, $\forall i\in\mathcal{I}$; the third ``$=$'' is obtained using $\rho_i=\sin{\alpha_i}$; the forth ``$=$'' is obtained using the formula $2\sin{a}\sin{(a+2b)}=\cos{2b}-\cos{(2a+2b)}$,
where $a\coloneqq\alpha_j$ and $b\coloneqq\sum_{k=1}^{j-1}{\alpha_k}$; the fifth ``$=$'' is obtained by cancellation; the last ``$=$'' is obtained using $2\sum_{i=1}^n{\alpha_i}=0~\text{mod}~2\pi$.
\\
Finally, $g_2(\xi)=0$ can be verified in a similar fashion; the details are omitted. 
This completes the sufficiency.
\\
\emph{Necessity:} When the rotation is counterclockwise, both $\alpha_i$ and $\alpha_i+\beta_i$ lie in the range $(0,\pi)$. As illustrated in Fig.~\ref{fig-sum-pi}, $\angle\mathcal{O}(i+1)i=\angle\mathcal{O}i(i+1)$, which follows $\pi/2-\alpha_i=(\alpha_i+\beta_i)-\pi/2$. Thus, we have $\alpha_i+(\alpha_i+\beta_i)=\pi$ for all $i$;~\eqref{eq-sum-pi-mod-2pi} is obtained. Combining 
\eqref{eq-trigonometric-relation} and~\eqref{eq-metric-relation} yields $2\sin{\alpha_i}/\rho_i=1/r$ for all $i$ with $r$ being the radius;~\eqref{eq-same-ratio} is obtained. 
The proof for clockwise rotation is similar and thus omitted. This completes the proof.
\endproof

\subsection{Proof of Proposition~\ref{thm-alpha-sum}} \label{app-alpha-sum}
The analysis when $n=2$ is simplified in the sense that $\alpha_2=\alpha_1+\beta_1$ and $\alpha_1=\alpha_2+\beta_2$. By~Proposition~\ref{thm-exist.-condt.}, we have $\alpha_1+\alpha_2=\alpha_1+(\alpha_1+\beta_1)=\pi$ for counterclockwise rotation and $\alpha_1+\alpha_2=\alpha_1+(\alpha_1+\beta_1)=-\pi$ for clockwise rotation. \\
Consider $n=3$ and suppose that the rotation is counterclockwise. As shown in~Fig.~\ref{fig-geomet.-config.-n=3}, there are two possible arrangements. We show that $\sum_i{\alpha_i}$ is equal to $\pi$ for the first arrangement and $2\pi$ for the second. As shown in~Fig.~\ref{fig-geomet.-config.-3cases}, there are three possible cases for the first arrangement. \\
(i) The center (of the circle) lies inside $\triangle{123}$. Then, the sum of the interior angles of $\triangle{123}$ can be computed as $\sum_{i=1}^3{2\left(\pi/2-\alpha_i\right)}=\pi$, which follows $\sum_{i=1}^3{\alpha_i}=\pi$. \\
(ii) The center lies on one side of $\triangle{123}$. Then, it is easy to verify that $\sum_{i=1}^3{\alpha_i}=\pi$ is still true. It does not matter which side the center lies on, as we can relabel the vehicles. \\
(iii) The center lies outside $\triangle{123}$. It still follows that $\sum_{i=1}^3{2\left(\pi/2-\alpha_i\right)}=\pi$ since $\pi/2-\alpha_3\in(-\pi/2,0)$ is clockwise and $2\left(\pi/2-\alpha_3\right)$ compensates $\sum_{i=1}^2{\left(\pi/2-\alpha_i\right)}$ when computing the sum of the interior angles of $\triangle{123}$. \\
Therefore, $\sum_{i=1}^3{\alpha_i}=\pi$ holds for all three cases. Then, if we reverse the direction of each vector $\mathtt{v}_i$ and view from inside the paper outward, the corresponding three cases of the second arrangement can be obtained. Denote the newly obtained bearing as $\alpha_i^\prime=\pi-\alpha_i$, then we have $\sum_{i=1}^3{\alpha_i^\prime}=\sum_{i=1}^3\left({\pi-\alpha_i}\right)=3\pi-\sum_{i=1}^3{\alpha_i}=2\pi$. \\
We now proceed with the case $n > 3$. Regard each triple $(1,i,i+1)$ as a virtual cyclic pursuit, where $i\in\{2,\ldots,n-1\}$, then $\sum_{i=1}^n{\alpha_i}$ can be obtained by summing up the corresponding sums of the bearings of these triples. Denote $\alpha_{1,j}$ as the virtual bearing from $1$'s heading to the heading that would take it directly towards $j$, and vice versa for $\alpha_{j,1}$, where $j\in\{3,\ldots,n-1\}$. \\
Suppose that the rotation is counterclockwise. We have shown that the sum of the bearings of a triple is $\pi$ or $2\pi$. Hence, the summation of the bearing sums of these triples is a positive multiple of $\pi$ and reads
\begin{equation*}\begin{aligned}
    &\left(\alpha_1+\alpha_2+\alpha_{3,1}\right) +\sum_{j=3}^{n-2}{\left(\alpha_{1,j}+\alpha_j+\alpha_{j+1,1}\right)}  \\
    &+\left(\alpha_{1,n-1}+\alpha_{n-1}+\alpha_{n}\right) \\
    &= \sum_{i=1}^n{\alpha_i}+\sum_{j=3}^{n-1}{\left(\alpha_{1,j}+\alpha_{j,1}\right)} = \sum_{i=1}^n{\alpha_i} + (n-3)\pi = \mathtt{q}\pi,
\end{aligned}\end{equation*}
where $\mathtt{q}\in\mathbb{Z}_+$, and $\alpha_{1,j}+\alpha_{j,1}=\pi$ follows from the previous analysis when $n=2$. \\
There are $n-2$ triples, and the value of $\mathtt{q}$ depends on the number of triples with a bearing sum of $\pi$ or equivalently $2\pi$. If all triples have a bearing sum of $\pi$, then $\mathtt{q}$ reaches the minimum $n-2$, which follows $\sum_{i=1}^n{\alpha_i}=\pi$. If all triples have a bearing sum of $2\pi$, then $\mathtt{q}$ reaches the maximum $2(n-2)$, which follows $\sum_{i=1}^n{\alpha_i}=(n-1)\pi$. Therefore, $\sum_{i=1}^n{\alpha_i}$ is a positive multiple of $\pi$ ranging from $\pi$ to $(n-1)\pi$. We thus verify \eqref{eq-sum-multiple} for counterclockwise rotation. The proof for clockwise rotation is analogous and thus omitted. This completes the proof.
\endproof

\begin{figure}
    \centering
    {
        \begin{tikzpicture}[scale=0.75]
            \draw[->, thick] (0,1.5) -- (-.75,1.5) node[left] {\scalebox{0.75}{$\mathtt{v}_1$}};
            \draw[->, thick] (-1.5,0) -- (-1.5,-.75)  node[below] {\scalebox{0.75}{$\mathtt{v}_2$}};
            \draw[->, thick] ({1.5/sqrt(2)},{-1.5/sqrt(2)}) -- ({1.5/sqrt(2)+.75/sqrt(2)}, {-1.5/sqrt(2)+.75/sqrt(2)}) node[right] {\scalebox{0.75}{$\mathtt{v}_3$}};
            \node at (0,0) [right] {\scalebox{0.75}{$z_\mathcal{O}$}};
        
            \draw[->, thick] (-.3,1.5) arc[start angle=-180, end angle=-135, radius=0.3];
            \node at (-.55, 1.25) {\scalebox{0.75}{$\alpha_1$}};
            \draw[->, thick] (-1.5,-.3) arc[start angle=-90, end angle=-22.5, radius=0.3];
            \node at (-1.2, -.5) {\scalebox{0.75}{$\alpha_2$}};
            \draw[->, thick] ({1.5/sqrt(2)+.3/sqrt(2)}, {-1.5/sqrt(2)+.3/sqrt(2)}) arc[start angle=45, end angle=112.5, radius=0.3];
            \node at ({1.5/sqrt(2)+.1/sqrt(2)}, {-1.5/sqrt(2)+.8/sqrt(2)}) {\scalebox{0.75}{$\alpha_3$}};
            
            \draw[thick] (0,1.5) -- (-1.5,0);
            \draw[thick] (-1.5,0) -- ({1.5/sqrt(2)},{-1.5/sqrt(2)});
            \draw[thick] ({1.5/sqrt(2)},{-1.5/sqrt(2)}) -- (0,1.5);
        
            \draw[dashed, thick] (0,0) -- (0,1.5);
            \draw[dashed, thick] (0,0) -- (-1.5,0);
            \draw[dashed, thick] (0,0) -- ({1.5/sqrt(2)},{-1.5/sqrt(2)});
        \end{tikzpicture}
    }
    {
        \begin{tikzpicture}[scale=0.75]
            \draw[->, thick] (0,1.5) -- (-.75,1.5) node[left] {\scalebox{0.75}{$\mathtt{v}_1$}};
            \draw[->, thick] (-1.5,0) -- (-1.5,-.75)  node[below] {\scalebox{0.75}{$\mathtt{v}_2$}};
            \draw[->, thick] (0,-1.5) -- (.75,-1.5) node[right] {\scalebox{0.75}{$\mathtt{v}_3$}};
            \node at (0,0) [right] {\scalebox{0.75}{$z_\mathcal{O}$}};
        
            \draw[->, thick] (-.3,1.5) arc[start angle=-180, end angle=-135, radius=0.3];
            \node at (-.55, 1.25) {\scalebox{0.75}{$\alpha_1$}};
            \draw[->, thick] (-1.5,-.3) arc[start angle=-90, end angle=-45, radius=0.3];
            \node at (-1.25, -.5) {\scalebox{0.75}{$\alpha_2$}};
            \draw[thick] (.25,-1.5) -- (.25,-1.25);
            \draw[->, thick] (.25,-1.25) -- (0,-1.25);
            \node at (.4, -1.1) {\scalebox{0.75}{$\alpha_3$}};
            
            \draw[thick] (0,1.5) -- (-1.5,0);
            \draw[thick] (-1.5,0) -- (0,-1.5);
            \draw[thick] (0,-1.5) -- (0,1.5);
        
            \draw[dashed, thick] (0,0) -- (0,1.5);
            \draw[dashed, thick] (0,0) -- (-1.5,0);
            \draw[dashed, thick] (0,0) -- (0,-1.5);
        \end{tikzpicture}
    }
    {
        \begin{tikzpicture}[scale=0.75]
            \draw[->, thick] (0,1.5) -- (-.75,1.5) node[left] {\scalebox{0.75}{$\mathtt{v}_1$}};
            \draw[->, thick] (-1.5,0) -- (-1.5,-.75)  node[below] {\scalebox{0.75}{$\mathtt{v}_2$}};
            \draw[->, thick] ({-1.5/sqrt(2)},{-1.5/sqrt(2)}) -- ({-1.5/sqrt(2)+.6/sqrt(2)}, {-1.5/sqrt(2)-.6/sqrt(2)}) node[right] {\scalebox{0.75}{$\mathtt{v}_3$}};
            \node at (0,0) [right] {\scalebox{0.75}{$z_\mathcal{O}$}};
        
            \draw[->, thick] (-.3,1.5) arc[start angle=-180, end angle=-135, radius=0.3];
            \node at (-.55, 1.25) {\scalebox{0.75}{$\alpha_1$}};
            \draw[->, thick] (-1.5,-.5) arc[start angle=-90, end angle=-67.5, radius=0.5];
            \node at (-1.05, -.45) {\scalebox{0.75}{$\alpha_2$}};
            \draw[->, thick] ({-1.5/sqrt(2)+.19/sqrt(2)}, {-1.5/sqrt(2)-.19/sqrt(2)}) arc[start angle=-45, end angle=67.5, radius=.19];
            \node at ({-1.5/sqrt(2)+.45}, {-1.5/sqrt(2)}) {\scalebox{0.75}{$\alpha_3$}};
    
            \draw[thick] (0,1.5) -- (-1.5,0);
            \draw[thick] (-1.5,0) -- ({-1.5/sqrt(2)},{-1.5/sqrt(2)});
            \draw[thick] ({-1.5/sqrt(2)},{-1.5/sqrt(2)}) -- (0,1.5);
        
            \draw[dashed, thick] (0,0) -- (0,1.5);
            \draw[dashed, thick] (0,0) -- (-1.5,0);
            \draw[dashed, thick] (0,0) -- ({-1.5/sqrt(2)},{-1.5/sqrt(2)});
        \end{tikzpicture}
    }
    \caption{Three scenarios of the first arrangement in Fig.~\ref{fig-geomet.-config.-n=3}.}
    \label{fig-geomet.-config.-3cases}
\end{figure}

\subsection{Proof of Proposition~\ref{thm-cot-sum}} \label{app-cot-sum}
We first study the case $\mathtt{p}=1$, i.e., $\sum_{i=1}^n{\alpha_i}=\pi$. Since $\alpha_i\in(0,\pi)$ for all $i$, there is at most one $j$ such that $\alpha_j\geq\pi/2$. Thus, $\cot{\alpha_j}\leq0$ and $\cot{\alpha_i}\geq0$ for all $i\neq j$. Since $\sum_{i=1}^n{\alpha_i}=\pi$, $\alpha_i<\pi-\alpha_j\leq\pi/2$ for all $i\neq j$, we have $\cot{\alpha_i}>\cot{(\pi-\alpha_j)}\Rightarrow\cot{\alpha_i}+\cot{\alpha_j}>0$  for all $i\neq j$. Thus, $\sum_{i=1}^n{\cot{\alpha_i}}>0$. 
\\
Consider the case $\mathtt{p}=n-1$, i.e., $\sum_{i=1}^n{\alpha_i}=(n-1)\pi$. When $n=3$, there are only two scenarios: (i) $\alpha_i\in(\pi/2,\pi)$ for all $i\in\{1,2,3\}$; (ii) Two angles lie in the range $(\pi/2,\pi)$, and the last lie in the range $(0,\pi/2]$. It is obvious that $\sum_{i=1}^3{\cot{\alpha_i}}<0$ for (i). For (ii), suppose $\alpha_1\leq\pi/2<\alpha_2\leq\alpha_3$, then it holds that $\cot{\alpha_3}\leq\cot{\frac{\alpha_2+\alpha_3}{2}} =\cot{\frac{2\pi-\alpha_1}{2}} <-\cot{\alpha_1}$,
which follows that $\cot{\alpha_1}+\cot{\alpha_3}<0$. Thus, $\sum_{i=1}^3{\cot{\alpha_i}}<0$ since $\cot{\alpha_2}<0$ as well. When $n>3$, by using the same idea in the proof of~Proposition~\ref{thm-alpha-sum}, we obtain
\begin{equation*}\begin{aligned}
    &\sum_{i=1}^n{\cot{\alpha_i}} 
    = \sum_{i=1}^n{\cot{\alpha_i}} + \sum_{j=3}^{n-1}{\left(\cot{\alpha_{1,j}}+\cot{\alpha_{j,1}}\right)} \\
    &=\left(\cot{\alpha_1}+\cot{\alpha_2}+\cot{\alpha_{3,1}}\right) \\
    &\quad +\sum_{j=3}^{n-1}{\left(\cot{\alpha_{1,j}}+\cot{\alpha_{j}}+\cot{\alpha_{j+1,1}}\right)} \\
    &\quad +\left(\cot{\alpha_{1,n-1}}+\cot{\alpha_{n-1}}+\cot{\alpha_{n}}\right),
\end{aligned}\end{equation*}
where the first ``='' holds since  $\alpha_{1,j}+\alpha_{j,1}=\pi$. We regard each triple, $(1,i,i+1)$ with $i\in\{2,3\ldots,n-1\}$, as a virtual cyclic pursuit. $\sum_{i=1}^n{\alpha_i}=(n-1)\pi$ implies that all triples have a bearing sum of $2\pi$; therefore, all triples have a negative value of $\sum{\cot(\cdot)}$, which guarantees $\sum_{i=1}^n{\cot{\alpha_i}}<0$. 
\\
Finally, when $\mathtt{p} \in (1, n-1) \subset \mathbb{Z}_+$, triples with a bearing sum of $\pi$ yield a positive value of $\sum{\cot(\cdot)}$, while those with a bearing sum of $2\pi$ yield a negative one. Consequently, the range of $\sum_{i=1}^n{\cot{\alpha_i}}$ is $\mathbb{R}$.
This completes the proof.
\endproof

\subsection{Proof of Lemma~\ref{thm-f-iinvariance}}\label{app-f-iinvariance}
\begin{figure*}
\begin{align}
    &\frac{\partial g_1(\xi)}{\partial\xi_i}f(\xi_i,\xi_{i+1}) = -v\left(\cos{\alpha_i}+\cos{(\alpha_i+\beta_i)}\right)\sin{(\alpha_i+\gamma_i)} + v\left( \sin{(\alpha_i+\beta_i)}-\sin{\alpha_i}\right)\cos{(\alpha_i+\gamma_i)} \label{eq-partial-g1f} \\
    & -\frac{2v}{\rho_i}\sin{\alpha_i}\sum_{j=i+1}^n{\rho_j\cos{(\alpha_j+\gamma_j)}}  +\frac{2v}{\rho_{i+1}}\sin{\alpha_{i+1}}\sum_{j=i+1}^n{\rho_j\cos{(\alpha_j+\gamma_j)}} -k\left(\cos{\alpha_i}+\cos{(\alpha_i+\beta_i)}\right)\cos{(\alpha_i+\gamma_i)} \nonumber \\ 
    & - \frac{k}{\rho_i}\left(\cos{\alpha_i}+\cos{(\alpha_i+\beta_i)}\right)\sum_{j=i+1}^n{\rho_j\cos{(\alpha_j+\gamma_j)}} +\frac{k}{\rho_{i+1}}\left(\cos{\alpha_{i+1}}+\cos{(\alpha_{i+1}+\beta_{i+1})}\right)\sum_{j=i+1}^n{\rho_j\cos{(\alpha_j+\gamma_j)}}. \nonumber \\
    &-\sum_{i=1}^n{\frac{2v}{\rho_i}\sin{\alpha_i}\sum_{j=i+1}^n{\rho_j\cos{(\alpha_j+\gamma_j)}}} +\sum_{i=1}^n{\frac{2v}{\rho_{i+1}}\sin{\alpha_{i+1}}\sum_{j=i+1}^n{\rho_j\cos{(\alpha_j+\gamma_j)}}} \label{eq-partial-g1f-comp1} \\
    &= - \frac{2v}{\rho_1}\sin{\alpha_1}\left(g_2(\xi)-\rho_1\cos{(\alpha_1+\gamma_1)}\right) - \sum_{i=2}^n{\frac{2v}{\rho_i}\sin{\alpha_i}\sum_{j=i+1}^n{\rho_j\cos{(\alpha_j+\gamma_j)}}} + \sum_{i=2}^n{\frac{2v}{\rho_i}\sin{\alpha_i}\sum_{j=i}^n{\rho_j\cos{(\alpha_j+\gamma_j)}}} \nonumber \\ &\quad +\frac{2v}{\rho_{n+1}}\sin{\alpha_{n+1}}\sum_{j=n+1}^n{\rho_j\cos{(\alpha_j+\gamma_j)}} = 
    - \frac{2v}{\rho_1}\sin{\alpha_1}g_2(\xi) + 2v\sum_{i=1}^n{\sin{\alpha_i}\cos{(\alpha_i+\gamma_i)}}. \nonumber \\
    & \frac{\partial g_2(\xi)}{\partial\xi_i}f(\xi_i,\xi_{i+1}) = -v\left(\cos{\alpha_i}+\cos{(\alpha_i+\beta_i)}\right)\cos{(\alpha_i+\gamma_i)} - v\left( \sin{(\alpha_i+\beta_i)}-\sin{\alpha_i}\right)\sin{(\alpha_i+\gamma_i)} \label{eq-partial-g2f}  \\
    & +\frac{2v}{\rho_i}\sin{\alpha_i}\sum_{j=i+1}^n{\rho_j\sin{(\alpha_j+\gamma_j)}} -\frac{2v}{\rho_{i+1}}\sin{\alpha_{i+1}}\sum_{j=i+1}^n{\rho_j\sin{(\alpha_j+\gamma_j)}} + k\left(\cos{\alpha_i}+\cos{(\alpha_i+\beta_i)}\right)\sin{(\alpha_i+\gamma_i)} \nonumber \\
    & + \frac{k}{\rho_i}\left(\cos{\alpha_i}+\cos{(\alpha_i+\beta_i)}\right)\sum_{j=i+1}^n{\rho_j\sin{(\alpha_j+\gamma_j)}} - \frac{k}{\rho_{i+1}}\left(\cos{\alpha_{i+1}}+\cos{(\alpha_{i+1}+\beta_{i+1})}\right)\sum_{j=i+1}^n{\rho_j\sin{(\alpha_j+\gamma_j)}}. \nonumber 
\end{align}
\hrulefill
\end{figure*}
The proof proceeds along the lines of \cite{marshall2005coordinated}. 
Rewrite $\frac{\partial g(\xi)}{\partial\xi}\hat{f}(\xi)$ as 
$\textstyle \sum_{i=1}^n\frac{\partial g(\xi)}{\partial\xi_i}f(\xi_i,\xi_{i+1})$. 
Denote $\gamma_i\coloneqq(i-1)\pi-\sum_{j=1}^{i-1}{\beta_j}$. Then $g_1(\xi)$ is rewritten as $g_1(\xi) = \sum_{i=1}^n{\rho_i\sin{(\alpha_i+\gamma_i)}}$.
The partial derivative ${\partial g_1(\xi)}\over{\partial\xi_i}$ reads 
\begin{align*}
    \left[\sin{(\alpha_i+\gamma_i)},~ \rho_i\cos{(\alpha_i+\gamma_i)},~ -\sum_{j=i+1}^n{\rho_j\cos{(\alpha_j+\gamma_j)}}\right].
\end{align*}
Notice that $\sum_{j={n+1}}^n{\rho_j\cos{(\alpha_j+\gamma_j)}}=0$. 
Multiplying $\frac{\partial g_1(\xi)}{\partial\xi_i}$ by $f(\xi_i,\xi_{i+1})$ gives~\eqref{eq-partial-g1f}. By summing the third and fourth terms of~\eqref{eq-partial-g1f} for indices $i$, we obtain~\eqref{eq-partial-g1f-comp1}. Then, summing the last term of~\eqref{eq-partial-g1f-comp1} and the first two terms of~\eqref{eq-partial-g1f} for indices $i$ yields
\begin{equation}\label{eq-partial-g1f-comp2}\begin{aligned}
    & v\sum_{i=1}^n{\left( \sin{(\alpha_i+\beta_i)}+\sin{\alpha_i}\right)\cos{(\alpha_i+\gamma_i)}} \\
    & -v\sum_{i=1}^n{\left(\cos{\alpha_i}+\cos{(\alpha_i+\beta_i)}\right)\sin{(\alpha_i+\gamma_i)}} \\
    & = v\sum_{i=1}^n{\left(\sin{(\beta_i-\gamma_i)}-\sin{\gamma_i}\right)} =v\sin{(\sum_{i=1}^n{\beta_i}-(n-1)\pi)} \\ 
    &= v\sin{(g_3(\xi)-\pi)} = -v\sin{(g_3(\xi))}. 
\end{aligned}\end{equation}
In a manner similar to~\eqref{eq-partial-g1f-comp1}, by summing the last two terms of~\eqref{eq-partial-g1f} for indices $i$ and noting that 
\begin{align*}
    \sum_{j=2}^n{\rho_j\cos{(\alpha_j+\gamma_j)}} = g_2(\xi)-
    \rho_1\cos{(\alpha_1+\gamma_1)},
\end{align*}
the summation can be concluded as
\begin{equation}\label{eq-partial-g1f-comp3}\begin{aligned}
    &-\frac{k}{\rho_1}\left(\cos{\alpha_1}+\cos{(\alpha_1+\beta_1)}\right)g_2(\xi) \\
    &+ k\sum_{i=1}^n{\left(\cos{\alpha_i}+\cos{(\alpha_i+\beta_i)}\right)\cos{(\alpha_i+\gamma_i)}}. 
\end{aligned}\end{equation}
Summing $n$ of the fifth term of~\eqref{eq-partial-g1f} cancels the last term of the above summation, and we finally obtain
\begin{equation*}\begin{aligned}
    \frac{\partial g_1(\xi)}{\partial\xi}\hat{f}(\xi) =-\omega_1g_2(\xi) -v\sin{(g_3(\xi))} = 0.
\end{aligned}\end{equation*}
The same rationale applies to $g_2(\xi)$, which can be rewritten as $g_2(\xi) = \sum_{i=1}^n{\rho_i\cos{(\alpha_i+\gamma_i)}}$. The partial derivative $\frac{\partial g_2(\xi)}{\partial\xi_i}$ reads
\begin{align*}
    \left[\cos{(\alpha_i+\gamma_i)},~ -\rho_i\sin{(\alpha_i+\gamma_i)},~ \sum_{j=i+1}^n{\rho_j\sin{(\alpha_j+\gamma_j)}}\right].
\end{align*}
Multiplying it by $f(\xi_i,\xi_{i+1})$ gives~\eqref{eq-partial-g2f}. As in~\eqref{eq-partial-g1f-comp1}, summing the third and fourth terms for indices $i$ gives
\begin{align*}
    \frac{2v}{\rho_1}\sin{\alpha_1}g_1(\xi) - 2v\sum_{i=1}^n{\sin{\alpha_i}\sin{(\alpha_i+\gamma_i)}}.
\end{align*}
By summing the last term of the above summation and the first two terms of~\eqref{eq-partial-g2f} for indices $i$, and proceeding as in \eqref{eq-partial-g1f-comp2}, we have
\begin{equation*}\begin{aligned}
    &-v\sum_{i=1}^n{\left(\cos{\alpha_i}+\cos{(\alpha_i+\beta_i)}\right)\cos{(\alpha_i+\gamma_i)}} \\
    &-v\sum_{i=1}^n{\left( \sin{(\alpha_i+\beta_i)}+\sin{\alpha_i}\right)\sin{(\alpha_i+\gamma_i)}} \\
    &= -v\sum_{i=1}^n{\left(\cos{(\beta_i-\gamma_i)}+\cos{\gamma_i} \right)} 
    = v\left(\cos{(g_3(\xi))}-1\right).
\end{aligned}\end{equation*}
As in~\eqref{eq-partial-g1f-comp3}, summing the last two terms of~\eqref{eq-partial-g2f} for indices $i$ yields
\begin{align*}
    &\frac{k}{\rho_1}\left(\cos{\alpha_1}+\cos{(\alpha_1+\beta_1)}\right)g_1(\xi) \\
    &-k\sum_{i=1}^n{\left(\cos{\alpha_i}+\cos{(\alpha_i+\beta_i)}\right)\sin{(\alpha_i+\gamma_i)}}.
\end{align*}
Finally, summing the fifth term of~\eqref{eq-partial-g2f} for indices $i$ cancels the last term of the above summation, and we end up with
\begin{equation*}\begin{aligned}
    \frac{\partial g_2(\xi)}{\partial\xi}\hat{f}(\xi) = \omega_1g_1(\xi) + v\left(\cos{(g_3(\xi))}-1\right) = 0.
\end{aligned}\end{equation*}
For $g_3(\xi)$, we can directly obtain 
\begin{align*}
    \frac{\partial g_3(\xi)}{\partial\xi}\hat{f}(\xi) = \sum_{i=1}^n{\dot{\beta}_i}=\sum_{i=1}^n{\omega_i-\omega_{i+1}} = 0
\end{align*}
since the indices $i+1$ are taken modulo $n$. 
\\
In conclusion,
$\frac{\partial g(\xi)}{\partial\xi}\hat{f}(\xi) = 0,~\forall\xi\in\mathcal{M}$.
It then follows as a consequence of Lemma~7 in~\cite{marshall2004formations} that the submanifold $\mathcal{M}$ is invariant under $\hat{f}$.
\endproof

\subsection{Proof of Lemma~\ref{thm-A-hat-pmjomega}}\label{app-A-hat-pmjomega}
Let $\varphi=\Phi(\xi)$ be the change of coordinates:
\begin{equation*}\begin{aligned}
    &\varphi_1=\rho_1,~\varphi_2=\alpha_1,~\ldots,~\varphi_{3n-3}=\beta_{n-1}, \\
    &\varphi_{3n-2}=g_1(\xi),~\varphi_{3n-1}=g_2(\xi),~\varphi_{3n}=g_3(\xi).
\end{aligned}\end{equation*}
Structure the new coordinates as $\varphi=[\varphi_{\mathrm{I}}^\top,\varphi_{\mathrm{II}}^\top]^\top$, where $\varphi_{\mathrm{I}}=[\varphi_1,\varphi_2,\ldots,\varphi_{3n-3}]^\top$ and $\varphi_{\mathrm{II}}=[\varphi_{3n-2},\varphi_{3n-1},\varphi_{3n}]^\top$. The coordinates in $\varphi_{\mathrm{II}}$ are the functions that define $\mathcal{M}$. In the new coordinates, it follows that
\begin{equation*}\begin{aligned}
    \dot{\varphi}_{\mathrm{I}} &= \left.\begin{bmatrix}
        I_{3n-3} & 0_{(3n-3)\times3}
    \end{bmatrix}\hat{f}(\xi)\right|_{\xi = \Phi^{-1}(\varphi)}, \\
    \dot{\varphi}_{\mathrm{II}} &= \begin{bmatrix}
        \displaystyle \frac{\partial g(\xi)}{\partial\xi}\hat{f}(\xi)
    \end{bmatrix}_{\xi = \Phi^{-1}(\varphi)}.
\end{aligned}\end{equation*}
The equilibrium $\bar{\varphi}=\Phi(\bar{\xi})$ is equal to $\bar{\xi}$, except that the last three components are instead zero. By computing the linearization at $\bar{\varphi}$, we get
\begin{equation*}\begin{aligned}
    \dot{\varphi}_{\mathrm{I}} &= \left.
        \begin{bmatrix}
        I_{3n-3} & 0_{(3n-3)\times3}
    \end{bmatrix}\begin{bmatrix}
        \displaystyle   \frac{\partial\hat{f}(\xi)}{\partial\varphi}
    \end{bmatrix}_{\xi=\Phi^{-1}(\varphi)}
    \right|_{\bar{\varphi}}\varphi \\
    &\coloneqq\begin{bmatrix}
        \hat{A}_{T_{\bar{\xi}}\mathcal{M}} & * 
    \end{bmatrix}\varphi, \\
    \dot{\varphi}_{\mathrm{II}} &= \left.
        \displaystyle \frac{\partial}{\partial\varphi}\begin{bmatrix}
            \displaystyle \frac{\partial g(\xi)}{\partial\xi}\hat{f}(\xi)
        \end{bmatrix}_{\xi=\Phi^{-1}(\varphi)}
    \right|_{\bar{\varphi}}\varphi \\
    &= \left.
        \displaystyle \frac{\partial}{\partial\varphi}
        \begin{bmatrix}
        \displaystyle    -\omega_1\varphi_{3n-1}  -v\sin{\varphi_{3n}} \\
        \displaystyle    \omega_1\varphi_{3n-2}  +v\left(\cos{\varphi_{3n}}-1\right) \\ 0
        \end{bmatrix}
    \right|_{\bar{\varphi}}\varphi \\
    &= \left[\begin{array}{ccc|ccc}
    0 & \ldots & 0 & 0 & -\bar{\omega} & -v \\
    0 & \ldots & 0 & \bar{\omega} & 0 & 0 \\
    0 & \ldots & 0 & 0 & 0 & 0 
    \end{array}\right]\varphi 
    = \begin{bmatrix}
        0_{3\times(3n-3)} & \hat{A}_{T_{\bar{\xi}}\mathcal{M}}^\star
    \end{bmatrix}\varphi.
\end{aligned}\end{equation*}
Readers are referred to the proof of~Lemma~\ref{thm-f-iinvariance} for the derivation of $\frac{\partial g(\xi)}{\partial\xi}\hat{f}(\xi)$. 
It is then clear that $\hat{A}_{T_{\bar{\xi}}\mathcal{M}}^\star$ has eigenvalues $\lambda_{1,2,3}=\{0,\pm j\bar{\omega}\}$.
\endproof

\subsection{Proof of Theorem~\ref{thm-stable-criterion}}\label{app-set-stable}
For all $\bar{\xi}\in\mathcal{C}$, since $\bar{s}\neq0$, it follows that
\begin{equation*}
    \mathcal{S}_\mathsf{s}(\bar{\xi
    }) \coloneqq \mspan{\left\{\left[\bar{\rho}_1,0,0,\ldots,\bar{\rho}_n,0,0\right]^\top\right\}}\subset\ker{\hat{A}}.
\end{equation*}
The corresponding zero eigenvalue stems from the system's degrees of freedom on $\mathcal{C}$, allowing circular formations of arbitrary radius, since $\bar{\xi}+\tilde{\xi}\in\mathcal{C},~\forall \tilde{\xi}\in\mathcal{S}_\mathsf{s}(\bar{\xi})$. \\
Further, write $\mathsf{v}_i^\mathsf{d} =\left[\bar{\rho}_i\cot{\bar{\alpha}_i}, 1, -2\right]^\top$. It holds for each pair of $A_i$ and $B_i$ that $A_i\mathsf{v}_i^\mathsf{d} = 0$ and $B_i\mathsf{v}_i^\mathsf{d}=0$,
and thus the following subspace is contained in $\ker{\hat{A}}$:
\begin{equation*}
    \mathcal{S}_{\mathsf{d},i}(\bar{\xi
    }) \coloneqq \mspan{\left\{\left[0_{3\times1}^\top,\ldots, {\mathsf{v}_i^\mathsf{d}}^\top,\ldots,0_{3\times1}^\top\right]^\top\right\}}\subset\ker{\hat{A}},
\end{equation*}
where $\mathsf{v}_i^\mathsf{d}$ is the $i$th $3$-tuple and the others are all zero. Each $\xi_i^\mathsf{d} \coloneqq \left[0_{3\times1}^\top,\ldots, {\mathsf{v}_i^\mathsf{d}}^\top,\ldots,0_{3\times1}^\top\right]^\top$
stems from the system's degrees of freedom on $\mathcal{C}$, allowing circular formations with arbitrary agent spacing. This intuition is justified using limiting analysis.
Let $\tilde{\xi}_i=[\tilde{\rho}_i,\tilde{\alpha}_i,\tilde{\beta}_i]^\top$ be a displacement that makes unicycle $i$ drifts along the orbit. Thus,  $2(\bar{\alpha}_i+\tilde{\alpha}_i)+(\bar{\beta}_i+\tilde{\beta}_i) = 2\bar{\alpha}_i+\bar{\beta}_i=\pi~\text{mod}~2\pi$, which follows that $\tilde{\alpha}_i\colon\tilde{\beta}_i=1\colon-2$. It also holds that ${\sin{(\bar{\alpha}_i+\tilde{\alpha}_i)}}/(\bar{\rho}_i+\tilde{\rho}_i)={\sin{\bar{\alpha}_i}}/{\bar{\rho}_i}$. Let $\tilde{\rho}_i\colon\tilde{\alpha}_i=\varepsilon\colon1$, it yields that $({\sin{(\bar{\alpha}_i+\tilde{\alpha}_i)}-\sin{\bar{\alpha}_i}})/{\tilde{\alpha}_i}=\varepsilon\cdot\left({\sin{\bar{\alpha}_i}}/{\bar{\rho}_i}\right)$, By taking the limit of the left-hand side as $\tilde{\alpha}_i\rightarrow0$, we obtain that $\varepsilon=\bar{\rho}_i\cot{\bar{\alpha}_i}$. \\
Recall that $\hat{A}$ can be transformed into~\eqref{eq-A-upper-form}, where $\hat{A}_{T_{\bar{\xi}}\mathcal{M}}^\star$ has solely imaginary eigenvalues $\lambda_{1,2,3} = \{0,\pm j2v\bar{s}\}$ for all $\bar{\xi}\in\mathcal{C}$; see Lemma~\ref{thm-A-hat-pmjomega}. Hence, $\hat{A}$ has at least $n+3$ solely imaginary eigenvalues. On the one hand, the submanifold corresponding to circular formation in relative coordinates has dimension $n$ only; see Remark~\ref{rmk-dimension-dismatch}. This apparent dimension mismatch is due to the use of redundant coordinates (say) $\rho_n$, $\alpha_n$, and $\beta_n$. With such redundancy, $\hat{A}$ can be derived directly with a cyclically interconnected structure relatively easier to analyze. Otherwise, we would have to derive an explicit expression of $\hat{A}_{T_{\bar{\xi}}\mathcal{M}}$ in~\eqref{eq-A-upper-form}, which is cumbersome and less tractable. On the other hand, the submanifold corresponding to circular formation in global coordinates has exactly dimension $n+3$; see Proposition~\ref{thm-exist.-condt.-global}. The three missing dimensions in the relative-coordinate representation are indeed the translation and rotation of the formation  in global coordinates. \\
Due to the zero eigenvalues, individual equilibria are no longer asymptotically stable, and our attention turns to the invariant set $\mathcal{C}$. By Propositions~\ref{thm-arrangements}--\ref{thm-alpha-sum}, the full invariant set $\mathcal{C}$ decomposes into a finite number of maximal connected subsets. Let $\mathcal{C}^*\subset\mathcal{C}$ be one such subset and $\bar{\xi}\in\mathcal{C}^*$ be a point of interest. Define $\mathcal{U}_\delta(\bar{\xi})\coloneqq\left\{\xi\in\mathcal{M}: \|\xi-\bar{\xi}\|<\delta\right\}$. If we choose
\begin{equation*}
    \delta < \inf_{\mathcal{C}^\prime\subset\mathcal{C} \setminus\mathcal{C}^*}{\dist{\left(\mathcal{C}^\prime,\mathcal{C}\right)}},
\end{equation*}
then $\mathcal{U}_\delta(\bar{\xi})\cap\left(\mathcal{C} \setminus\mathcal{C}^*\right)=\emptyset$.
Similarly, denote $\mathcal{U}_\epsilon(\bar{\xi})\coloneqq\left\{\bar{\xi}^\prime\in\mathcal{C} : \|\bar{\xi}^\prime-\bar{\xi}\|<\epsilon\right\}$, $\mathcal{U}_\epsilon(\bar{\xi})\cap\left(\mathcal{C} \setminus\mathcal{C}^*\right)=\emptyset$ if $\epsilon$ is chosen in the same way as $\delta$. \\
We have shown that $\hat{A}$ has at least $n+3$ solely imaginary eigenvalues for all $\bar{\xi}^\prime\in\mathcal{C}$, and the submanifold corresponding to circular formation has exactly dimension $n+3$ in global coordinates. Thus, if the remaining $2n-3$ eigenvalues of $\hat{A}$ at $\bar{\xi}\in\mathcal{C}^*$ have negative real parts, then
$\mathcal{U}_\epsilon(\bar{\xi})$ becomes a center manifold, and there exist sufficiently small $\epsilon>0$ and $\delta>0$, consistent with the above choices, such that
\begin{equation*}
    \xi(0)\in \mathcal{U}_\delta(\bar{\xi})
    \Rightarrow \lim_{t\rightarrow\infty}{\dist{\left(\xi(t),\mathcal{U}_\epsilon(\bar{\xi})\right)}} = 0.
\end{equation*}
That is, points outside $\mathcal{C}$ and sufficiently close to $\bar{\xi}\in\mathcal{C}^*$ will converge to a neighborhood within $\mathcal{C}^*$. \\
It is preferable to have uniform stability on $\mathcal{C}^*$. Since $\mathcal{C}^*$ is maximal connected, $\mathcal{C}^*=\bigcup_{\bar{\xi}\in\mathcal{C}^*}{\mathcal{U}_\epsilon(\bar{\xi})}$. Hence, if the spectral condition is satisfied by every $\bar{\xi}\in\mathcal{C}^*$, then $\mathcal{C}^*$ becomes a center manifold and is locally asymptotically stable. This completes the proof.
\endproof

\subsection{Proof of Lemma~\ref{thm-A-reduced}}\label{app-A-reduced}
$\bar{s}\neq0$ can be factored out of $\hat{A}$, so that some entries become functions of $\cot{\bar{\alpha}_i}$ (Note that $\bar{c}_i/\bar{s}=\bar{c}_i\bar{\rho}_i/\sin{\bar{\alpha}_i}=\cot{\bar{\alpha}_i}$). The explicit dependence on $\bar{\rho}_i$'s of the other entries is not appealing and could be removed through algebraic operations that might not affect the stability of $\hat{A}$. The underlying reason is that, intuitively, the stability is dependent on the bearing spacing $\bar{\alpha}_i$'s, not on the radius of the circle. It will be shown that such an algebraic operation is scaling each coordinate $\rho_i$ by $1/\bar{\rho}_i$. We now introduce a new coordinate $\gamma_i\coloneqq2\alpha_i+\beta_i$ to replace $\beta_i$. The above procedures are equivalent to applying transformation $\hat{A}_U = U\hat{A}U^{-1}$ with $U=\blkdiag{\left(U_1,U_2,\ldots,U_n\right)}$ and each $U_i$ given as 
\begin{align*}
    U_i= \begin{bmatrix}
        1/{\bar{\rho}_i} & 0 & 0 \\ 0 & 1 & 0 \\ 0 & 2 & 1
    \end{bmatrix}.
\end{align*}
In the matrix $\hat{A}_U$, it holds for all $i$ that the $(3i-1)$th column is $\cot{\bar{\alpha}_i}$ times the $(3i-2)$th column. Thus, consider another transformation $\hat{A}_V= V\hat{A}_UV^{-1}$ with $V=\blkdiag{\left(V_1,V_2,\ldots,V_n\right)}$, where 
\begin{align*}
    V_i = \begin{bmatrix}
        1 & -\cot{\alpha}_i & 0 \\ 0 & 1 & 0 \\
        0 & 0 & 1
    \end{bmatrix}.
\end{align*}
In the matrix $\hat{A}_V$, it holds for all $i$ that the $(3i-1)$th column is null. Hence, it can be transformed by permutations into block triangular form~\eqref{eq-A-lower-form}, and $\hat{A}_R$ is then obtained as in~\eqref{eq-A-reduced}.
\endproof

\subsection{Proof of Theorem~\ref{thm-eigen.-poly.}}\label{app-eigen.-poly.}
Given that $k\bar{s}<0\Rightarrow k\cdot\sign{(\bar{s})}=-2v<0$ and  $v=1$, $\bar{A}_i$ and $\bar{B}$ can be rewritten as
\begin{align*}
    \bar{A}_i& = |\bar{s}|\begin{bmatrix}
        0 & \sign{(\bar{s})}z_i^2 \\
        -2\cdot\sign{(\bar{s})} & -2z_i
    \end{bmatrix}, \\
    \bar{B} &= |\bar{s}|\begin{bmatrix}
        0 & 0 \\ 2\cdot\sign{(\bar{s})} & -2
    \end{bmatrix},
\end{align*}
where $z_i \coloneqq 1+\sign{(\bar{s})}\cot{\bar{\alpha}_i}$. Denote $N=\hat{A}_R/|\bar{s}|$ and write $M\coloneqq\lambda I-N$. Then $M$ has the form
\begin{align*}
    M &=\left[\begin{array}{cccc:c}
        Z_1 & Y & & &\\
        & Z_2 & Y & & \\
        & & \ddots & \ddots & \\
        & & & Z_{n-1} & Y \\
        \cdashline{1-5}
        Y & & & & Z_n
        \end{array}\right]= \begin{bmatrix}
        A & B \\ C & D\end{bmatrix}, \quad\text{where} \\
    Z_i &= \begin{bmatrix}
        \lambda & -\sign{(\bar{s})}z_i^2 \\
        2\cdot\sign{(\bar{s})} & \lambda+2z_i
    \end{bmatrix}, ~
    Y = \begin{bmatrix}
        0 & 0 \\ -2\cdot\sign{(\bar{s})} & 2
    \end{bmatrix}.
\end{align*}
It is well known that
\begin{align*}
    \det{M} = \det{A}\det{\left(D-CA^{-1}B\right)}.
\end{align*}
Since $A$ is block upper-triangular, we have
\begin{align*}
    \det{A} &= \prod_{i=1}^{n-1}{\det{Z_i}} = \prod_{i=1}^{n-1}{\left(\lambda^2+2z_i\lambda+2z_i^2\right)},  \\
    A^{-1} &= \begin{bmatrix}
    Z_1^{-1} & *_{12} & \ldots & *_{1(n-1)} \\
    & Z_2^{-1} & \ldots & *_{2(n-1)} \\
    & & \ddots & \vdots \\
    & & & Z_{n-1}^{-1}
\end{bmatrix}, \quad\text{where} \\
    *_{ij} &\coloneqq (-1)^{j-i}Z_i^{-1}\prod_{k=i+1}^j{YZ_{k}^{-1}},~1\leq i<j\leq n-1, \\
    &\prod_{k=i+1}^j{YZ_{k}^{-1}} \coloneqq YZ_{i+1}^{-1}YZ_{i+2}^{-1}\ldots YZ_j^{-1}.
\end{align*}
Note that $Z_i$ is now $2\times2$ instead of $3\times3$ prior to model reduction. Otherwise, computing its inverse and carrying out the subsequent derivation would be more cumbersome. Further, we have
\begin{align*}
    CA^{-1}B &= Y*_{1(n-1)}Y 
    = (-1)^{n}\left(\prod_{i=1}^{n-1}{YZ_i^{-1}}\right)Y, \quad\text{where} \\
    YZ_i^{-1} &= \frac{1}{\det{Z_i}}\begin{bmatrix}
    0 & 0 \\ *_i & 2(\lambda-z_i^2) \end{bmatrix}, \\
    *_i &\coloneqq-2\cdot\sign{(\bar{s})}(\lambda+2z_i+2).
\end{align*}
It then follows that
\begin{equation*}\begin{aligned}
    &\left(\prod_{i=1}^{n-1}{YZ_i^{-1}}\right)Y \\
    &= \frac{1}{\prod_{i=1}^{n-1}{\det{Z_i}}}\left(\prod_{i=1}^{n-1}\begin{bmatrix}
        0 & 0 \\ *_i & 2(\lambda-z_i^2)
    \end{bmatrix}\right)\begin{bmatrix}
        0 & 0 \\ -2\cdot\sign{(\bar{s})} & 2
    \end{bmatrix} \\
    &=\frac{1}{\det{A}}\begin{bmatrix}
        0 & 0 \\ * & 2^{n-1}\prod_{i=1}^{n-1}{\left(\lambda-z_i^2\right)}
    \end{bmatrix}\begin{bmatrix}
        0 & 0 \\ -2\cdot\sign{(\bar{s})} & 2
    \end{bmatrix} \\
    &= \frac{1}{\det{A}}\cdot{2^n\prod_{i=1}^{n-1}{\left(\lambda-z_i^2\right)}}\begin{bmatrix}
        0 & 0 \\ -\sign{(\bar{s})} & 1
    \end{bmatrix}.
\end{aligned}\end{equation*}
No need to specify $*$ as it plays no role in the computation. We are now ready to compute $\det{\left(D-CA^{-1}B\right)}$:
\begin{equation*}\begin{aligned}
    & \det{\left(D-CA^{-1}B\right)} = \det{\left(Z_n-(-1)^{n}\left(\prod_{i=1}^{n-1}{YZ_i^{-1}}\right)Y\right)} \\
    &= \begin{vmatrix}
        \lambda & -\sign{(\bar{s})}z_n^2 \\
        2\cdot\sign{(\bar{s})}+ \sign{(\bar{s})}/\det{A}\cdot* & \lambda+2z_n - 1/\det{A}\cdot* 
    \end{vmatrix} \\
    &= \lambda^2+2z_n\lambda+2z_n^2-\left(\lambda-z_n^2\right)/{\det{A}}\cdot* \\ 
    &= \lambda^2+2z_n\lambda+2z_n^2-{1}/{\det{A}}\cdot\prod_{i=1}^{n}{\left(2z_i^2-2\lambda\right)},
    \label{eq-det-D-CinvAB}
\end{aligned}\end{equation*}
where $*\coloneqq(-2)^{n}\prod_{i=1}^{n-1}{\left(\lambda-z_i^2\right)}$.
Combining the above results, we arrive at the following polynomial:
\begin{equation*}\begin{aligned}\label{eq-implicit-det.}
    \det{\left(\lambda I-N\right)}  = \prod_{i=1}^n{\left(\lambda^2+2z_i\lambda+2z_i^2\right)} - \prod_{i=1}^{n}{\left(2z_i^2-2\lambda\right)}.
\end{aligned}\end{equation*}
The matrix $\hat{A}_R$ has eigenvalues $\{0,\pm j2\bar{s}\}$, which correspond to $\{0,\pm j2\}$ for $N=\hat{A}_R/|\bar{s}|$. The rest of the theorem then follows.
\\
Finally, while~Lemma~\ref{thm-A-hat-pmjomega} confirms that $\{0,\pm j2\}$ are the roots of~$P_{\bar{\xi}}(\lambda)$, we provide a direct proof. It is evident that zero is a root. Then we show that $j2$ is also a root. It follows that
\begin{align*}
    P_{\bar{\xi}}(\lambda) = 0 ~\Leftrightarrow~\prod_{i=1}^n{\frac{\lambda^2+2z_i\lambda+2z_i^2}{2z_i^2-2\lambda}}=1.
\end{align*}
Suppose $\lambda=j2$. It yields that
\begin{align*}
    \left\|\lambda^2+2z_i\lambda+2z_i^2\right\|^2 &= \left(2z_i^2-4\right)^2+\left(4z_i\right)^2 = 4z_i^4+16, \\
    \left\|2z_i^2-2\lambda\right\|^2 &= \left(2z_i^2\right)^2 + \left(-4\right)^2 = 4z_i^4+16.
\end{align*}
Therefore, we have 
\begin{align*} 
    \left|\frac{2z_i^2-4+j4z_i}{2z_i^2-j4}\right|=1.
\end{align*}
Here, $|\cdot|$ denotes the complex modulus. Thus, it holds for some $\phi_i$ that
\begin{equation*}\begin{aligned}
    e^{j\phi_i} &= \frac{2z_i^2-4+j4z_i}{2z_i^2-j4} = \frac{1}{z_i^4+4}\left(z_i^2-2+j2z_i\right)\left(z_i^2+j2\right) \\
    &= \frac{1}{z_i^4+4}\left(\left(z_i^2-2\right)z_i^2-4z_i+j\left(2(z_i^2-2)+2z_i^3\right)\right).
\end{aligned}\end{equation*}
Then we can obtain that
\begin{equation*}\begin{aligned}
    \cot{\phi_i} &= \frac{\left(z_i^2-2\right)z_i^2-4z_i}{2(z_i^2-2)+2z_i^3} \\
    &= \frac{z_i\left(z_i-2\right)\left(z_i^2+2z_i+2\right)}{2\left(z_i-1\right)\left(z_i^2+2z_i+2\right)} = \frac{z_i^2-2z_i}{2z_i-2}.
\end{aligned}\end{equation*}
On the other hand, recall that $z_i\coloneqq 1+\sign{(\bar{s})}\cot{\bar{\alpha}_i}=1+\cot{|\bar{\alpha}_i|}$ with $|\bar{\alpha}_i|\in(0,\pi)$; see~Remark~\ref{rmk-cot|alpha|}, and we have
\begin{equation*}\begin{aligned}
    \cot{2|\bar{\alpha}_i|} &= \frac{\cot^2{|\bar{\alpha_i}|-1}}{2\cot{|\bar{\alpha}_i|}} \\
    &= \frac{\left(z_i-1\right)^2-1}{2\left(z_i-1\right)} = \frac{z_i^2-2z_i}{2z_i-2}.
\end{aligned}\end{equation*}
Therefore, we can conclude that $\phi_i=2|\bar{\alpha}_i|$.  
By~Proposition~\ref{thm-alpha-sum}, we have $\sum_{i=1}^n{\phi_i}=2\sum_{i=1}^{n}{|\bar{\alpha}_i|}=2\texttt{p}\pi$, where $\texttt{p}\in[1,n-1]\subset\mathbb{Z}_+$. Then we can derive 
\begin{equation*}\begin{aligned}
    \prod_{i=1}^n{\frac{2z_i^2-4+j4z_i}{2z_i^2-j4}} =\prod_{i=1}^n{e^{j\phi_i}}
    =e^{j2\texttt{p}\pi}=1,
\end{aligned}\end{equation*}
which confirms that $j2$ is a root of $P_{\bar{\xi}}(\lambda)$. The proof for $-j2$ is analogous and thus omitted. 
\endproof

\subsection{Proof of Theorem~\ref{thm-stable-n=3}}\label{app-stable-n=3}
Let $v=1$ without loss of generality. The polynomial~\eqref{eq-implicit-poly.} when $n=3$ has the form
\begin{equation*}\begin{aligned}\label{eq-poly.-n=3}
    P_{\bar{\xi}}(\lambda)
    &= \lambda^6+b_5\lambda^5+b_4\lambda^4+b_3\lambda^3 + b_2\lambda^2+b_1\lambda \\
    &= \left(\lambda^3+4\lambda\right)\left(\lambda^3+a_2\lambda^2+a_1\lambda+a_0\right), 
\end{aligned}\end{equation*}
with the coefficients $a$'s and $b$'s yet to be determined. \\
To obtain the coefficients $a$'s, we expand~\eqref{eq-implicit-poly.} and derive the expression of the coefficients $b$'s:
\begin{equation*}\begin{aligned}
    b_5 & = 2\left(z_1+z_2+z_3\right),~
    b_4 = 2\left(z_1+z_2+z_3\right)^2, \\
    b_3 & = 4z_1^2\left(z_2+z_3\right) + 4z_2^2\left(z_1+z_3\right) + 4z_3^2\left(z_1+z_2\right) \\
    &~\quad +8\left(z_1z_2z_3+1\right), \\
    b_1 & = 8z_1z_2z_3\left(z_1z_2+z_2z_3+z_1z_3\right) +8\left(z_1^2z_2^2+z_2^2z_3^2+z_1^2z_3^2\right).
\end{aligned}\end{equation*}
We omit $b_2$ since it will not be used. Then, the coefficients $a$'s can be obtained as
\begin{equation*}
    a_2 = b_5,~
    a_1 = b_4-4,~ a_0 =b_3-4b_5 = b_1/4.
\end{equation*}
Recall that $z_i\coloneqq 1+\sign{(\bar{s})}\cot{\bar{\alpha}_i}=1+\cot{|\bar{\alpha}_i|}$ with $|\bar{\alpha}_i|\in(0,\pi)$; see~Remark~\ref{rmk-cot|alpha|}. In other words, we can focus solely on counterclockwise rotation. Then we show that the coefficients a's satisfy the Routh-Hurwitz criterion for all $\bar{\xi}\in\mathcal{C}_\pi^-=\left\{\bar{\xi}^\prime\in\mathcal{C}^- : \sum_{i=1}^3{|\bar{\alpha}_i|=\pi}\right\}$. \\
Since $\sum_{i=1}^3{|\bar{\alpha}_i|}=\pi\Rightarrow\sum_{i=1}^{3}{\cot{|\bar{\alpha}_i}|}>0$, it yields that $a_2 = 2\left(z_1+z_2+z_3\right) = 6+2\sum_{i=1}^{3}{\cot{|\bar{\alpha}_i}|}>0$,
and $ a_1 = 2\left(z_1+z_2+z_3\right)^2-4 = 2\left({a_2}/{2}\right)^2 - 4 > 0$.
Next, it follows that
\begin{equation}\begin{aligned}
    &a_0= 2z_1z_2z_3\left(z_1z_2+z_2z_3+z_1z_3\right) \\
    &~\quad\quad + 2\left(z_1^2z_2^2+z_2^2z_3^2+z_1^2z_3^2\right). \label{eq-a0-1} 
\end{aligned}\end{equation}
Since $\cot{\left(|\bar{\alpha}_1|+|\bar{\alpha}_2|\right)} = \cot{\left(\pi-|\bar{\alpha}_3|\right)} = -\cot{|\bar{\alpha}_3|}$, we have $\cot{|\bar{\alpha}_1|}\cot{|\bar{\alpha}_2|}+\cot{|\bar{\alpha}_2|}\cot{|\bar{\alpha}_3|} +\cot{|\bar{\alpha}_1|}\cot{|\bar{\alpha}_3|}=1$. Hence,
\begin{equation}\begin{aligned}
    &z_1z_2+z_2z_3+z_1z_3 \\ &= \cot{|\bar{\alpha}_1|}\cot{|\bar{\alpha}_2|}+\cot{|\bar{\alpha}_2|}\cot{|\bar{\alpha}_3|} +\cot{|\bar{\alpha}_1|}\cot{|\bar{\alpha}_3|}  \\ &\quad+2\left(\cot{|\bar{\alpha}_1|}+\cot{|\bar{\alpha}_2|}+\cot{|\bar{\alpha}_3|}\right) + 3 \\ 
    &=2\left(\cot{|\bar{\alpha}_1|}+\cot{|\bar{\alpha}_2|}+\cot{|\bar{\alpha}_3|}\right)+4 \\
    &=2\left(z_1+z_2+z_3\right)-2. \label{eq-a0-2}
\end{aligned}\end{equation}
Since $\sum_{i=1}^3{|\bar{\alpha}_i|}=\pi$ and $|\bar{\alpha}_i|\in (0,\pi)$, there can be at most one $i$ such that $|\bar{\alpha}_i|\geq\pi/2\Rightarrow\cot{|\bar{\alpha}_i|}\leq0$. Furthermore, if such $|\bar{\alpha}_i|>3\pi/4$, then $z_i=1+\cot{|\bar{\alpha}_i|}<0$. Suppose $0<|\bar{\alpha}_1|\leq|\bar{\alpha}_2|<\pi/4<3\pi/4<|\bar{\alpha}_3|$, which yields $z_1,z_2>0$ and $z_3<0$. Then,  combining~\eqref{eq-a0-1}--\eqref{eq-a0-2} yields $a_0=2\left(z_1z_2+z_2z_3+z_1z_3\right)^2-4z_1z_2z_3>0$.
On the other hand, if there is no $i$ such that $|\bar{\alpha}_i|>3\pi/4$, then all $z_i\geq0$, with at most one being zero and at least two greater than one. For this case, we can directly conclude that $a_0>0$ based on~\eqref{eq-a0-1}. \\
It remains to show $a_2a_1-a_0>0$, and we have
\begin{equation}\begin{aligned}
    &a_2a_1-a_0 = b_5\left(b_4-4\right)-\left(b_3-4b_5\right) =b_5b_4-b_3 \\
    &= 4\left(z_1+z_2+z_3\right)^3 -4z_1^2\left(z_2+z_3\right) - 4z_2^2\left(z_1+z_3\right) \\
    &\quad - 4z_3^2\left(z_1+z_2\right) -8\left(z_1z_2z_3+1\right) \\
    &= 4\left(z_1^3+z_2^3+z_3^3\right) +8z_1^2\left(z_2+z_3\right) + 8z_2^2\left(z_1+z_3\right) \\
    &\quad + 8z_3^2\left(z_1+z_2\right) +16z_1z_2z_3-8. \label{eq-a2a1-a0}
\end{aligned}\end{equation}
Again, if there is no $i$ such that $|\bar{\alpha}_i|>3\pi/4$, then all $z_i\geq0$, with at most one being zero and at least two greater than one. For this case, it is not difficult to prove $a_2a_1-a_0>0$ by verifying that $4\left(z_1^3+z_2^3+z_3^3\right)-8>0$; the details are omitted. Otherwise, suppose $0<|\bar{\alpha}_1|\leq|\bar{\alpha}_2|<\pi/4<3\pi/4<|\bar{\alpha}_3|$. Expanding~\eqref{eq-a2a1-a0} yields
\begin{equation*}\begin{aligned}\label{eq-a2a1-a0-expand}
    &a_2a_1-a_0 = 4\sum_{i=1}^3{\left(1+3\cot{|\bar{\alpha}_i|}+3\cot^2{|\bar{\alpha}_i|}+\cot^3{|\bar{\alpha}_i|}\right)} \\
    &\quad +48\sum_{i=1}^3{\cot{|\bar{\alpha}_i|}}+16\sum_{i=1}^3{\cot^2{|\bar{\alpha}_i|}}+80 \\
    &\quad +8\sum_{i=1}^3{\cot^2{|\bar{\alpha}_i|}\sum_{j\neq i}{\cot{|\bar{\alpha}_j|}}} \\
    &\quad +16\left(2+\sum_{i=1}^3{\cot{|\bar{\alpha}_i|}}+\prod_{i=1}^3{\cot{|\bar{\alpha}_i}|}\right) -8. 
\end{aligned}\end{equation*}
It follows from $\cot{|\bar{\alpha}_1|}\geq\cot{|\bar{\alpha}_2|}>\left|\cot{|\bar{\alpha}_3|}\right|>1$ that $\sum_{i=1}^3{\cot^2{|\bar{\alpha}_i|}\sum_{j\neq i}{\cot{|\bar{\alpha}_j|}}}>0$ and
$\sum_{i=1}^3{\cot^q{|\bar{\alpha}_i|}}>0$, $q\in\{1,2,3\}$. On the other hand, $16\prod_{i=1}^3{\cot{|\bar{\alpha}_i|}}<0$ due to $\cot{|\bar{\alpha}_3|}<0$, but the following sum can be constructed to offset this negative term:
\begin{equation*}\begin{aligned}\label{eq-sum-cubic}
    &\frac{16}{6}\sum_{i=1}^3{\cot^3{|\bar{\alpha}_i|}}+\frac{16}{6}\times3\sum_{i=1}^3{\cot^2{|\bar{\alpha}_i|}\sum_{j\neq i}{\cot{|\bar{\alpha}_j|}}} \\
    &+\frac{16}{6}\times6\prod_{i=1}^3{\cot{|\bar{\alpha}_i|}} = \frac{16}{6}\left(\sum_{i=1}^3{\cot{|\bar{\alpha}_i|}}\right)^3>0.
\end{aligned}\end{equation*} 
Thus, $\tilde{P}_{\bar{\xi}}(\lambda)$ is Hurwitz; the remaining three eigenvalues of $\hat{A}$ have negative real parts. It then follows from Theorem~\ref{thm-stable-criterion} that $\mathcal{C}_\pi^-$ is locally asymptotically stable.
\endproof

\balance
\subsection{Proof of Theorem~\ref{thm-negeigen}}\label{app-negeigen}
While the idea behind the proof is simple, the detailed derivation is tedious; therefore, we provide only a proof sketch summarizing its main steps and rationale.
\\
Still, we assume $v=1$ without loss of generality. Then, since $|\bar{\alpha}_i|={\pi}/{n}$ for all $i$, the polynomial~\eqref{eq-implicit-poly.} is reduced to $P(\lambda) = \left(\lambda^2+2z\lambda+2z^2\right)^n-\left(2z^2-2\lambda\right)^n$,
where $z=1+\cot{(\pi/n)}$. The roots of $P(\lambda)$ are explicitly given by
\begin{equation}\label{formu-polyroot}
\lambda = -\left(z + \omega^{i-1}\right) \pm \sqrt{\left(z +\omega^{i-1}\right)^2-2z^2\left(1-\omega^{i-1}\right)},
\end{equation}
where $\omega^{i-1}\coloneqq e^{j2(i-1)\pi/n}\in\mathbb{C}$. Let $\phi_i\coloneqq2(i-1)\pi/n$, then $\omega^{i-1}=\cos{\phi_i}+j\sin{\phi_i}$. 
\\
For ease of reference later, we introduce additional notations. First, write
\begin{align}
    -\left(z+\omega^{i-1}\right) &= a_{i,1} + jb_{i,1}, \label{a_i1} \\
    \sqrt{\left(z+\omega^{i-1}\right)^2-2z^2\left(1-\omega^{i-1}\right)} &= a_{i,2}+jb_{i,2}.
\end{align}
Since $z=1+\cot{(\pi/n)}>1$ and $-1\leq\re{\omega^{i-1}}\leq1$, we have $a_{i,1}<0$ for all $i$. With a slight abuse of notation, let $z=a+jb$. It holds that $\re{\sqrt{z}}>0$ if $a>0$. Hence, we can conclude that $a_{i,2}>0$  for all $i$ due to
$\re{\left\{\left(z+\omega^{i-1}\right)^2-2z^2\left(1-\omega^{i-1}\right)\right\}}>0$. 
\\
Then we write
\begin{equation}\begin{aligned}
    \Phi_i &\coloneqq \left(z+\omega^{i-1}\right)^2-2z^2\left(1-\omega^{i-1}\right) \nonumber \\
    &= a_{i,2}^2-b_{i,2}^2 + j2a_{i,2}b_{i,2}. \label{Z_i} 
\end{aligned}\end{equation}
The next step is to verify that all solutions to~\eqref{formu-polyroot} lie in the open left half-plane, with the exception of one zero and a pair of $\pm j2$. We summarize the following four cases.
\\
(i) For $i=1$, one can easily yield $\lambda=0$ and $\lambda=-2\left(z+1\right)<0$ as the solutions to~\eqref{formu-polyroot}.
\\
(ii) For $i=2$, the idea is to verify that $\lambda = j2$ is one solution to~\eqref{formu-polyroot}. Suppose this is the case, then it follows that $a_{2,1}+a_{2,2}=0$ and $b_{2,1}+b_{2,2}=2$. Under these conditions, it suffices to prove that $a_{2,2}^2-b_{2,2}^2 = a_{2,1}^2 - \left(2-b_{2,1}\right)^2$ and $2a_{2,2}b_{2,2} = 2a_{2,1}\left(b_{2,1}-2\right)$.
The detailed (rather tedious) derivation is omitted. Thus, $\lambda=j2$ is one of the roots. The other root has a negative real part equal to $a_{2,1}-a_{2,2}=a_{2,1}-(-a_{2,1})=2a_{2,1}<0$.
\\
(iii) For \(i = n\), the analysis follows similarly to the last case. Here, the idea is to verify that $\lambda=-j2$ is one of the roots. 
Then, the other root has a negative real part equal to $2a_{n,1}<0$.
\\
(iv) For $i=3,4,\ldots,n-1$, write $\lambda = a_{i,1} + jb_{i,1} \pm \left( a_{i,2}+jb_{i,2} \right)$. We aim to prove that $a_{i,1}\pm a_{i,2}<0$, i.e., the roots have negative real parts. Since $a_{i,1}<0$ and $a_{i,2}>0$, it is equivalent to showing $a_{i,1}^2-a_{i,2}^2>0$.
\\
From~\eqref{a_i1} we have $a_{i,1}^2=\left(z+\cos{\phi_i}\right)^2$, where $\phi_i=2(i-1)\pi/n$. From~\eqref{Z_i} we can derive
\begin{equation*}
    a_{i,2}^2 = \frac{1}{2}\left( \re{\Phi_i} + \sqrt{ \left( \re{\Phi_i} \right)^2 + \left( \im{\Phi_i} \right)^2 } \right).
\end{equation*}
By expanding~\eqref{Z_i}, we can yield
\begin{subequations}\label{Z_i-re-im}\begin{align}
    \re{\Phi_i} &= -z^2 + 2\cos{\phi_i}\left(z^2+z\right) + \cos{2\phi_i}, \\
    \im{\Phi_i} &= 2\sin{\phi_i}\left(z^2+z\right) + \sin{2\phi_i}.
\end{align}\end{subequations}
Thus, we need to verify the following inequality
\begin{equation*}\begin{aligned}
    &2\left(z+\cos{\phi_i}\right)^2 - \re{\Phi_i} - \sqrt{ \left( \re{\Phi_i} \right)^2 + \left( \im{\Phi_i} \right)^2 } \\
    &= z^2\left(2-2\cos{\phi_i}\right) + z^2+2z\cos{\phi_i}+1 \\
    &\quad - \sqrt{ \left( \re{\Phi_i} \right)^2 + \left( \im{\Phi_i} \right)^2 } \\
    &= z^2\left(2-2\cos{\phi_i}\right) + \left(z+\cos{\phi_i}\right)^2 + \sin^2{\phi_i} \\
    &\quad - \sqrt{ \left( \re{\Phi_i} \right)^2 + \left( \im{\Phi_i} \right)^2 } > 0.
\end{aligned}\end{equation*}
We can easily tell that $z^2\left(2-2\cos{\phi_i}\right) + \left(z+\cos{\phi_i}\right)^2 + \sin^2{\phi_i}>0 \Rightarrow 2\left(z+\cos{\phi_i}\right)^2 - \re{\Phi_i}>0$.
Therefore, $a_{i,1}^2-a_{i,2}^2>0$ is equivalent to the condition $\Delta_{i} \coloneqq \left( 2\left(z+\cos{\phi_i}\right)^2 - \re{\Phi_i} \right)^2 - \left( \re{\Phi_i} \right)^2 - \left( \im{\Phi_i} \right)^2 >0$.
Substituting~\eqref{Z_i-re-im} into the condition yields a quadratic (in the form of $f(\cos{\phi_i})=a\cos^2{\phi_i}+b\cos{\phi_i}+c$). $\Delta_{i}>0$ can be verified by analyzing the vertex and discriminant of that quadratic. 
\\
Concluding all cases, all roots of $P(\lambda)$ (with the exception of one zero and a pair of $\pm j2$) have negative real parts; the remaining $2n-3$ eigenvalues of $\hat{A}$ have negative real parts. The statement then follows from Theorem~\ref{thm-stable-criterion}, together with arguments similar to those in its proof.
\endproof

\end{document}